\documentclass[a4paper,reqno]{amsart}

\RequirePackage{doi}
\usepackage{hyperref}

\usepackage{amsmath,amsthm,amssymb}
\usepackage{mathrsfs}
\usepackage[shortlabels]{enumitem}
\usepackage{graphicx}
\usepackage[font=small]{caption}
\usepackage{xcolor}
\usepackage{url}

\newcommand{\N}{\mathbb{N}}
\newcommand{\R}{{\mathbb{R}}}
\newcommand{\C}{{\mathbb{C}}}
\newcommand{\Z}{{\mathbb{Z}}}

\newcommand{\dd}{{{\rm d}}}
\newcommand{\ii}{{\rm i}}

\newcommand{\ov}{\overline}
\newcommand\wt{\widetilde}

	\newcommand{\la}{\lambda}
\newcommand{\La}{\Lambda}
\newcommand{\eps}{\varepsilon}

\newcommand{\spp}{\sigma_{\rm p}}

\newcommand{\Dom}{{\operatorname{Dom}}}
\newcommand{\Ker}{{\operatorname{Ker}}}
\newcommand{\Ran}{{\operatorname{Ran}}}

\newcommand{\Rank}{{\operatorname{rank}}}
\renewcommand{\Re}{\operatorname{Re}}
\renewcommand{\Im}{\operatorname{Im}}

\newcommand{\dist}{\operatorname{dist}}

\newcommand{\Tr}{\operatorname{Tr}}

\newcommand{\Num}{\operatorname{Num}}

\newcommand{\BigO}{\mathcal{O}}

\newcommand{\lspan}{{\operatorname{span}}}

\newcommand{\gs}{\gtrsim}

\theoremstyle{plain}

\newtheorem{theorem}{Theorem}[section]
\newtheorem{lemma}[theorem]{Lemma}

\newtheorem{proposition}[theorem]{Proposition}

\theoremstyle{definition}
\newtheorem{example}[theorem]{Example}
\newtheorem{remark}[theorem]{Remark}
\newtheorem{asm-sec}[theorem]{Assumption}

\newcommand\cB{\mathcal B}

\newcommand\cD{\mathcal D}

\newcommand\cG{\mathcal G}
\newcommand\cH{\mathcal H}

\newcommand\cM{\mathcal M}

\newcommand\cS{\mathcal S}

\numberwithin{equation}{section}

\begin{document}

\date{March 17, 2025}

\title[Local form subordination without a power decay]{Local form subordination without a power decay and a criterion of Riesz basesness}

\author{Boris Mityagin}

\address[Boris Mityagin]{Department of Mathematics, The Ohio State University, 231 West 18th Ave,
	Columbus, OH 43210, USA}
\email{mityagin.1@osu.edu, boris.mityagin@gmail.com}

\dedicatory{In memory of Nikolai Vasilevski}

\author{Petr Siegl}

\address[Petr Siegl]{Institute of Applied Mathematics, Graz University of Technology, Steyrergasse 30, 8010 Graz, Austria}

\email{siegl@tugraz.at}



\thanks{We gratefully acknowledge the support of BIRS for the Research in Teams stay 23rit099 in Banff, August 20 -- September 2, 2023, during which this work was initiated.}

\begin{abstract}
We revisit the local form subordination condition on the perturbation of a self-adjoint operators with compact resolvent, which is used to show the Riesz basis property of the eigensystem of the perturbed operator. Our new assumptions and new proof allow for establishing the Riesz basis property also in the case of slow and non-monotone decay in this subordination condition.  
\end{abstract}

\maketitle

\section{Introduction}
\label{sec:intro}

	Let $\cH$ be a separable Hilbert space and let $A$ be a self-adjoint operator with compact resolvent and simple eigenvalues $\{\mu_n\}_{n \in \N} \subset (0,\infty)$, ordered increasingly. Let $\{\psi_n\}_{n \in \N}$ be normalized eigenvectors of $A$ related to the eigenvalues $\mu_n$, 
\begin{equation}\label{A.def}
A \psi_n = \mu_n \psi_n, 
\quad 
\|\psi_n\| =1, \quad n \in \N = \{1,2,3,\dots\};
\end{equation}
notice that eigenvectors $\psi_n$ are in fact not determined uniquely by \eqref{A.def} as for any $s \in \R$, $e^{\ii s} \psi_n$, $n \in \N$, satisfy \eqref{A.def} as well. Nonetheless, for any admissible choice of $\{\psi_n\}_{n \in \N}$ in sense of \eqref{A.def},  $\{\psi_n\}_{n \in \N}$ forms an orthonormal basis of $\cH$.

This paper deals with the Riesz basis property of the eigensystem of non-symmetric perturbations of $A$;~i.e.~of 
\begin{equation}\label{T.def.intro}
T=A+V.
\end{equation}
The existing results are based on analysis of two ingredients: the size of the gaps of the unperturbed eigenvalues $\{\mu_n\}_{n \in \N}$ and the strength of the perturbation $V$ (typically not bounded). For instance, if the gaps grow, i.e.~$\mu_{n+1}-\mu_n \rightarrow + \infty$, and the perturbation $V$ is bounded, the result of Kato (\cite[Thm.~V.4.15a, Lem.~V.4.17a]{Kato-1966}, \cite{Kato-1967-73}) yields that the system of eigenvectors and possibly an additional finite number of root (or associated) vectors of $T$ contains a Riesz basis. Further results, relying on more precise information on the eigenvalues' gaps and including also relatively bounded and $p$-subordinated perturbations $V$ (in the operator or form sense), can be found e.g.~in \cite[Chap.~XIX.2]{DS3}, \cite[Chap.~I]{Markus-1988}, \cite{Clark-1968-25,Agranovich-1994-28,Wyss-2010-258}; see also \cite{Xu-2005-210, Zwart-2010-249,Motovilov-2017-8,Motovilov-2019-53} and \cite{Djakov-2010-283,Djakov-2012-263} where in particular perturbations of 1D Dirac operators (having constant eigenvalue gaps) are analyzed.

A recent progress initiated in \cite{Adduci-2012-10}, in particular in the case of constant gaps, is based on employing the Hilbert transform as an important technical tool and a more precise way to measure the strength of perturbation  via the effect of $V$ on the eigenvectors of $A$. Namely, it is assumed that $V$ satisfies 
\begin{equation}\label{B.loc.sub}
	\|V \psi_j \| \leq \eta_j, \quad j \in \N,
\end{equation}
or a form-version
\begin{equation}\label{B.loc.sub.form}
|\langle V \psi_j, \psi_k \rangle | \leq \omega_j \omega_k, \quad j,k \in \N,
\end{equation}
where the non-negative sequences $\{\eta_j\}_{j \in \N}$ or $\{\omega_j\}_{j \in \N}$ decay sufficiently fast relative to the eigenvalues' gaps; these conditions are also known as local (form) subordination; see \cite{Shkalikov-2010-269,Adduci-2012-73,Mityagin-2016-106,Shkalikov-2016-71,Mityagin-2019-139}. In particular, if the gaps behave as
\begin{equation}\label{gaps.reg}
\mu_{n+1} - \mu_n \gs n^{\gamma-1} 
\end{equation}
with some $\gamma>0$ and all sufficiently large $n \in \N$, then the eigensystem of $T$ contains a Riesz basis if 
\begin{equation}
\eta_j = o(j^{\gamma-1}), \quad j \to \infty,
\end{equation}
or, for the form version, if 
\begin{equation}\label{omega.reg}
\exists \alpha \in \R \text{ with } 2\alpha+\gamma>1: \quad  \omega_j = \BigO(j^{-\alpha}), \quad j \to \infty;	
\end{equation}
see \cite{Adduci-2012-10,Adduci-2012-73,Mityagin-2016-106,Mityagin-2019-139} for details.

The form local subordination \eqref{B.loc.sub.form} allows for perturbations that are merely sesquilinear forms $v : \Dom A^\frac 12 \times  \Dom A^\frac 12 \to \C$ or, in other words, operators from $(\Dom A^\frac 12,\|A^\frac12 \cdot\|)$ to its dual $(\Dom A^\frac 12,\|A^\frac12 \cdot\|)^*$. In applications to differential operators like the self-adjoint anharmonic oscillators
\begin{equation}
-\partial_x^2+|x|^\beta, \quad \beta  \geq 2, 
\end{equation}
in $L^2(\R)$, the Riesz basis property is established e.g.~for form perturbations generated by the potential $V = V_1 + V_2 + V_3 + V_4$ where
\begin{equation}\label{A.V.intro}
	\begin{aligned}
		&\exists \delta >0, \quad  |x|^{\frac{2-\beta}{2} +\delta} V_1(x) \in L^\infty(\R),
		\\ 
		&\exists p \in [1,\infty), \quad V_2 \in L^p(\R),
		\\
		&\exists s \in [0 ,\tfrac{\beta-1}{2 \beta}), \quad V_3  \in W^{-s,2}(\R),  
		\\
		&\exists \{\nu_k\}_{k \in \Z} \in \ell^1(\Z), \quad \exists \{x_k\}_{k \in \Z} \subset \R, \quad V_4  = \sum_{k \in \Z} \nu_k \, \delta(x-x_k);   
	\end{aligned}
\end{equation}
see~\cite{Mityagin-2016-106,Mityagin-2019-139} for more details and further examples.

In this paper we revisit the conditions on the gaps \eqref{gaps.reg} and the local form subordination~\eqref{B.loc.sub.form} with \eqref{omega.reg}. It is known that for the gap behavior \eqref{gaps.reg}, the decay condition \eqref{omega.reg} cannot be weakened to $\alpha$ with $2\alpha+\gamma=1$. In this case, the eigensystem of $T$ does not contain even a basis in general; see \cite[Sec.~5.1]{Mityagin-2019-139} inspired by the examples in \cite[Sec.~6.3]{Adduci-2012-10} and \cite[Sec.~8.1]{Adduci-2012-73}. On the other hand, the proofs in \cite{Mityagin-2016-106,Mityagin-2019-139} suggest that a slower decay than in \eqref{omega.reg} would be sufficient; nevertheless, the regular power-like behavior of the gaps and the sequence $\{\omega_j\}_{j \in \N}$ is used essentially in the given arguments. Our goal here is to improve such restrictions and find a new and weaker version of the condition~\eqref{omega.reg} for a general behavior of the gaps. 

We define the perturbed operator $T$ via the form sum; see Section~\ref{ssec:T.def}. The perturbation term $KVK$ in the resolvent of $T$, namely
\begin{equation}\label{res.repr.intro}
	(z-T)^{-1} = K(z)(I-K(z)VK(z))^{-1}K(z), \quad z \in \rho(A) \cap \rho(T),
\end{equation}
where $K(z)$ is a square root of $(z-A)^{-1}$, can be estimated as 
\begin{equation}\label{KVK.intro}
\|K(z)VK(z)\| \leq \|K(z)VK(z)\|_{\rm HS} \leq
\sum_{j =1}^\infty \frac{\omega_j^2}{|z-\mu_j|};
\end{equation}
see \eqref{K.def}, \eqref{Tz.res.dec} and \eqref{B.HS} below for details. Thus, as a minimal requirement on $\{\omega_j\}_{j \in \N}$, one could impose the summability condition 
\begin{equation}\label{asm.om.l1.intro}
\sum_{j=1}^\infty \frac{\omega_j^2}{\mu_j} < \infty.
\end{equation} 
Notice that \eqref{asm.om.l1.intro} is satisfied if \eqref{omega.reg} holds since \eqref{gaps.reg} implies that 
$\mu_j \gs j^\gamma$ for all sufficiently large $j$; see~\cite[Lemma~4.1]{Mityagin-2019-139}. Under \eqref{asm.om.l1.intro}, the resolvent of the perturbed operator $T$ is compact; see Proposition~\ref{prop:T.basic.def}. Moreover, if for some $p>0$
\begin{equation}\label{om.Sp.intro}
\sum_{j=1}^\infty \frac{1}{\mu_j^p}< \infty,
\end{equation}
i.e.~if $A^{-1}$ is in the Schatten class $\cS_p$, then the eigensystem of $T$ is complete; see Proposition~\ref{prop:T.complete}. Nevertheless, Example \ref{ex:counter} below for $\mu_j=j$, $j \in \N$, which is a yet another adjustment of \cite[Sec.~6.3]{Adduci-2012-10}, shows that \eqref{asm.om.l1.intro} alone is not sufficient for the Riesz basis property of the eigensystem of $T$ in general (not even for a basis).

Considering \eqref{KVK.intro} and the resolvent representation \eqref{res.repr.intro} again, another possible restriction on $\{\omega_j\}_{j \in \N}$ is that  
\begin{equation}\label{om.|HT|.intro.discs}
\sup_{|z-\mu_n| = r_n} \ \sum_{j =1}^\infty \frac{\omega_j^2}{|z-\mu_j|} = o(1), \quad n \to \infty,
\end{equation} 
where 
\begin{equation}
r_n := \frac12 \dist(\mu_n, \sigma(T) \setminus \{\mu_n\}), \quad n \in \N;
\end{equation}
i.e.~$r_n$ is the half of the distance of $\mu_n$ to the remaining eigenvalues of $T$. Indeed, under \eqref{om.|HT|.intro.discs}, all sufficiently large eigenvalues of $T$ remain simple and localized in a neighborhood of the eigenvalues of $A$; see Proposition~\ref{prop:localization}. In fact, instead of \eqref{om.|HT|.intro.discs}, we work with an equivalent condition 
\begin{equation}\label{om.|HT|.intro}
\cG[\{\omega_j^2\}_{j \in \N}](n):=\sum_{\substack{j =1 \\ j \neq n}}^\infty \frac{\omega_j^2}{|\mu_n-\mu_j|} + \frac{\omega_n^2}{r_n} = o(1), \quad n \to \infty
\end{equation} 
and, as our main result (Theorem~\ref{thm:RB.slow}), we establish the Riesz basis property of the eigensystem of $T$ if \eqref{B.loc.sub.form}, \eqref{asm.om.l1.intro}, \eqref{om.Sp.intro} and \eqref{om.|HT|.intro} hold.  We remark that, instead of \eqref{om.|HT|.intro}, it is also possible to assume only that, for a sufficiently small (fixed) $\eps>0$, there exists $N_\eps \in \N$, such that for all $n>N_\eps$, we have $|\cG[\{\omega_j^2\}_{j \in \N}](n)|<\eps$.

The condition \eqref{om.|HT|.intro} can be viewed as a transform $\cG$ of the sequence $\{\omega_j^2\}_{j \in \N}$ which for $\mu_j=j$, $j \in \N$, is close to the discrete Hilbert transform of $\{\omega_j^2\}_{j \in \N}$ (that is the sum in \eqref{om.|HT|.intro} without the absolute value in the denominator). The decay of the transformed sequence is satisfied in particular if \eqref{gaps.reg} and \eqref{omega.reg} hold, since for $2 \alpha + \gamma>1$
\begin{equation}\label{om.|HT|.intro.rates}
	\cG[\{j^{-2\alpha}\}_{j \in \N}](n) = 
	\begin{cases}
	\BigO(n^{-(2\alpha+\gamma-1)} \log n),  &\alpha \leq 1/2,
	\\[1mm]
	\BigO(n^{-\gamma}), & \alpha >1/2,
	\end{cases}
	\qquad n \to \infty;
\end{equation}  
see \cite[Lemma~4.1, 4.3]{Mityagin-2019-139}. 

Moreover, it follows from \eqref{om.|HT|.intro.rates} and \cite[Lemma~4.1, 4.3]{Mityagin-2019-139} that under \eqref{gaps.reg} and \eqref{omega.reg}, another transform of the sequence $\{\omega_j^2\}_{j \in \N}$, namely
\begin{equation}\label{G.transf.2.intro}
\wt \cG[\{\omega_j^2\}_{j \in \N}](k):=\sum_{\substack{n =1 \\ n \neq k}}^\infty \frac{r_n}{|\mu_k-\mu_n|} \sum_{\substack{j =1 \\ j \neq n}}^\infty \frac{\omega_j^2}{|\mu_n-\mu_j|}, \quad k \in \N,
\end{equation}
decays as $k \to \infty$; this property was essentially used in the previous proofs in \cite{Mityagin-2016-106,Mityagin-2019-139}. In fact, we show in Section~\ref{sec:fast} that if $\wt \cG[\{\omega_j^2\}_{j \in \N}]$ in \eqref{G.transf.2.intro} is a bounded sequence, a substantial simplification of the proof of Theorem~\ref{thm:RB.slow} is possible.
	
The main condition \eqref{om.|HT|.intro} is a rather implicit restriction on the sequence $\{\omega_j\}_{j \in \N}$. In Section~\ref{sec:ex}, we discuss the special case $\mu_j=j$, $j \in \N$, and find sufficient conditions both for the decay of  $\cG[\{\omega_j^2\}_{j \in \N}]$ in \eqref{om.|HT|.intro} and the boundedness of $\wt \cG[\{\omega_j^2\}_{j \in \N}]$ in \eqref{G.transf.2.intro} in this case. 

In Section~\ref{ssec:om.dec}, we show that for a \emph{monotone decreasing} $\{\omega_j\}_{j \in \N}$ the summability condition $\{\omega_j^2/j\}_{j \in \N} \in \ell^1(\N)$ (see \eqref{asm.om.l1.intro}) implies that \eqref{om.|HT|.intro} holds. Thus an example of such decreasing $\{\omega_j\}_{j \in \N}$ is 
\begin{equation}\label{om.lnln.intro}
\omega_j =\frac{1}{(\log j)^\frac12 (\log \log j)^{a} }, \quad j \geq 2, \ a > \frac 12;
\end{equation}
see Example~\ref{ex:lnln} for details. On the other hand, $\wt \cG[\{\omega_j^2\}_{j \in \N}]$ from \eqref{G.transf.2.intro} is bounded if 
\begin{equation}
\omega_j = \frac{1}{ (\log j)^a}, \qquad j \geq 2, \ a > 1;
\end{equation}
(and $\wt \cG[\{\omega_j^2\}_{j \in \N}]$ is not bounded for $a=1$); see Example~\eqref{ex:ln.fast} for details. These examples also illustrate the difference between requiring the decay of $\cG[\{\omega_j^2\}_{j \in \N}]$ in \eqref{om.|HT|.intro} and the boundedness of $\wt \cG[\{\omega_j^2\}_{j \in \N}]$ in \eqref{G.transf.2.intro}.

Finally, in Section~\ref{ssec:om.gaps}, we discuss \emph{non-monotone} sequences $\{\omega_j\}_{j \in \N}$ and show that the condition \eqref{om.|HT|.intro} is satisfied for $\{\omega_j\}_{j \in \N}$ with sufficiently large gaps in its support and an arbitrarily slow decay of its non-zero elements.

The paper is organized as follows. In Section~\ref{sec:prelim} we collect preliminaries. In Section~\ref{sec:result}, we introduce the perturbed operator $T$ and present our results. The proofs are given in Section~\ref{sec:proofs}. In Section~\ref{sec:fast} we consider the case of a faster decay of $\{\omega_j^2\}_{j \in \N}$. Finally, in Section~\ref{sec:ex}, we discuss the special case $\mu_j=j$, $j \in \N$. We find sufficient conditions for the main condition \eqref{om.|HT|.intro} and give examples of sequences satisfying the latter. We also present an example showing that the summability condition \eqref{asm.om.l1.intro} alone is not sufficient for the Riesz basis property.

\section{Preliminaries}
\label{sec:prelim}

We recall some standard notions and results. We follow mostly \cite{Markus-1988}, \cite{Kato-1966}, \cite[\S I.2]{Gohberg-1969} and \cite{DS2}. 


Let $T$ be a closed operator and let $\Lambda$ be a bounded closed subset of $\sigma(T)$ which is isolated (i.e.  $\sigma(T) \setminus \Lambda$ is closed). If $U$ is a bounded domain in $\C$ with a rectifiable boundary such that $\Lambda \subset U$ and $(\sigma(T) \setminus \Lambda) \cap \ov{U} = \emptyset$, let $P_{\Lambda}$ be the corresponding Riesz projections, i.e.
\begin{equation}
	P_{\Lambda} = \frac{1}{2 \pi \ii} \int_{\partial U} (z-T)^{-1} \, \dd z.
\end{equation}

A system of projections $\{Q_j\}_{j\in J}$ is called \emph{disjoint} if 
\begin{equation}\label{Q.disj.def}
	Q_j Q_k = \delta_{jk} Q_j, \quad j,k \in J.
\end{equation}
The disjointness of Riesz projections corresponding to disjoint subsets of spectrum follows by standard arguments; see e.g.~\cite[I.\S2.2]{Markus-1988}, \cite[Chap.~III.6.4 and 6.5]{Kato-1966} or \cite[\S I.2]{Gohberg-1969}. In detail, we have the following.
\begin{lemma}\label{lem:disj}
	If $\Lambda_j$, $j=1,2$, are bounded, closed, isolated and disjoint subsets of $\sigma(T)$, then the corresponding Riesz projections $P_{\Lambda_j}$$, j=1,2$, are disjoint.
\end{lemma}

Let $T$ be an operator with compact resolvent. Then any eigenvalue $\la \in \sigma(T)$ is isolated, for the corresponding Riesz projections we write $P_\la:=P_{\{\la\}}$, the corresponding root subspace $\Ran (P_{\lambda})$ is finite dimensional and it is equal to $\cup_{n=1}^\infty \Ker (A-\la)^n$. Any bounded subset $\La$ of $\sigma(T)$ is closed, isolated and $\Ran(P_\La)$ is the direct sum of the root subspaces corresponding to the eigenvalues in~$\La$.

A system of projections $\{Q_j\}_{j\in J}$ is said to be \emph{complete} if the following implication holds
\begin{equation}
\langle Q_j x,y \rangle =0 \quad \forall x\in \mathcal{H}, \ \forall j\in J    \quad \implies \quad y=0;	
\end{equation}
equivalently, it is complete if any  $x\in\mathcal H$ can be
approximated with an arbitrary accuracy by linear combinations of
$x_k \in \Ran (Q_k)$, $k\in J$.

For an operator with resolvent in a Schatten class, the completeness of the system of its Riesz projections follows if the resolvent norm can be controlled (in the sense of \eqref{res.pol.gr} below) on sufficiently many rays in $\C$ (depending on the Schatten class). For details and proof see \cite[Cor.XI.9.31]{DS2}.

\begin{theorem}\label{thm:compl}
	Let $T$ be a densely defined closed operator with $\rho(T) \neq \emptyset$ and for some $\la_0 \in \rho(T)$
	\begin{equation}
		(\la_0-T)^{-1} \in \cS_p
	\end{equation} 
	with $p \in (0,\infty)$. Let $\gamma_1, \dots, \gamma_s$ be non-overlapping differentiable arcs having a limit direction at infinity and suppose that no adjacent pair of arcs forms an angle as great as $\pi/p$ at infinity. Suppose that, for some $K \in \Z$, the resolvent of $T$ satisfies 
	\begin{equation}\label{res.pol.gr}
		\|(\la-T)^{-1}\| = \BigO(|\la|^K)
	\end{equation}
	as $\la \to \infty$ along each $\gamma_i$, $i=1,\dots, s$. Then the system of Riesz projections $P_{\lambda}$, $\la \in \sigma(T)$, is complete. 
\end{theorem}

Finally, our Theorem~\ref{thm:RB.slow} is based on the following criterion; for the proof see \cite[Chap.~6]{Gohberg-1969}, \cite[\S 6]{Shkalikov-2016-71} or \cite[Thm.~A]{Motovilov-2017-8}. 

\begin{theorem}\label{thm:RB.criteria}
	Let $\{Q_j\}_{j\in J}$ be a
	system of projections in a separable Hilbert space $\mathcal{H}$. The
	following statements are equivalent:
	\begin{enumerate}[\upshape i)]
		\item There exists a boundedly invertible operator $W \in \cB(\cH)$ and a complete system of disjoint orthogonal projections $\{Q_j^0\}_{j \in J}$ such that
		$$Q_j=W^{-1}Q_j^0 W,  \qquad j\in J.$$
		
		\item \label{item:RB} The system of projections $\{Q_j\}_{j\in J}$ is disjoint, complete and, for each $x \in \cH$
		\begin{equation}\label{0}
			\sum_{j\in J}|\langle Q_jx,x \rangle|<\infty.
		\end{equation}
	\end{enumerate}
\end{theorem}

\section{Results}
\label{sec:result}
In Section~\ref{ssec:T.def}, we introduce the perturbed operator $T$. In Section \ref{ssec:main.res}, we present our main results on the localization of eigenvalues and the Riesz basis property of the eigensystem. The proofs are given in Section~\ref{sec:proofs}.

\subsection{Perturbed operator $T$ and completeness of its eigensystem}
\label{ssec:T.def}

We assume that the operator $A$ is as in \eqref{A.def} with $\sigma(A)=\{\mu_n\}_{n\in \N} \subset (0,\infty)$. 

To define the perturbed operator $T$ we employ the representation theorems for closed sectorial forms and perturbation theory; see~\cite[Sec.VI]{Kato-1966}. By the second representation theorem (see~\cite[Thm.VI.2.23]{Kato-1966}), the operator $A$ is associated with the form
\begin{equation}\label{a.def}
	a[f]:=a(f,f) = \|A^\frac12 f\|^2, \quad \Dom(a)=\Dom(A^\frac12).
\end{equation}
The form domain equipped with its natural scalar product
\begin{equation}
	\cD_a := (\Dom(a), \langle A^\frac12 \cdot, A^\frac12 \cdot \rangle)	
\end{equation}
is a Hilbert space. 

Let $v$ be a sesquilinear form with $\Dom(v) = \cD_a$ such that 
\begin{equation}\label{asm:v.loc.sub}
	|v (\psi_j, \psi_k) | \leq \omega_j \omega_k, \quad j, k \in \N,
\end{equation}
with a sequence $\{\omega_j\}_{j \in \N} \subset [0,\infty)$ satisfying
\begin{equation}\label{asm:om.l1}
	\sum_{j=1}^\infty \frac{\omega_j^2}{ \mu_j} < \infty.
\end{equation}
\begin{proposition}\label{prop:T.basic.def}
	Let the operator $A$ be as in \eqref{A.def}, let $a$ be the associated form in \eqref{a.def}. Let $v$ be a sesquilinear form such that \eqref{asm:v.loc.sub} and \eqref{asm:om.l1} are satisfied. Then the form $v$ is bounded on $\cD_a$, moreover, for any $\eps>0$, there exists $C_\eps>0$ such that for all $f \in \cD_a$
	\begin{equation}\label{v.eps.a}
		|v[f]| \leq \eps (a[f] + C_\eps\|f\|^2),
	\end{equation}
	i.e.~$v$ is a relatively bounded perturbation of the form $a$ with the bound $0$. Thus the form
	\begin{equation}
		t := a + v, \quad \Dom(t) :=\Dom(a)
	\end{equation}
	is closed and sectorial and hence it defines, via the first representation theorem~\cite[Thm.~VI.2.1]{Kato-1966}, an m-sectorial operator $T$ with compact resolvent. 
	\item \label{prop:T.basic.ii}
\end{proposition}

We remark that, as a consequence of the boundedness of $v$ on $\cD_a$ under \eqref{asm:v.loc.sub} and \eqref{asm:om.l1}, any sesquilinear form $v$ defined initially only on $\lspan\{\psi_k\}_{k \in \N}$ for which \eqref{asm:v.loc.sub} holds with the sequence $\{\omega_j\}_{j \in \N}$ satisfying \eqref{asm:om.l1}, can be extended to a bounded form on $\cD_a$; see also the estimate \eqref{v.bdd} below. Moreover, one can view the form $v$ alternatively as a bounded operator $V : \cD_a \to \cD_a^*$ defined via $v(f,g) := \langle Vf,g \rangle_{\cD_a^* \times \cD_a}$ for all $f,g \in \cD_a$; this can be used to give a precise meaning to \eqref{B.loc.sub.form} in cases when $V$ is not an operator in $\cH$.

The operator $T$ can be also described as follows. We introduce the operators 
\begin{equation}\label{K.def}
	K(z) f := \sum_{k \in \N} (z-\mu_k)^{-\frac 12} \langle f, \psi_k \rangle \psi_k, \quad z \in \rho(A), \quad f \in \cH,
\end{equation}
where, for $0\neq w \in \C$ and $s\in \R$, the $s$-power of $w$ is taken as $w^s := |w|^s e^{\ii s \arg w}$ with $-\pi < \arg w \leq \pi$. Notice that
\begin{equation}
	\|K(z)\| = \max_{k \in \N}{ \frac{1}{|z-\mu_k|^\frac12}}, \quad z \in \rho(A), 
\end{equation}
and by \eqref{K.def}
\begin{equation}\label{K^2}
	K(z)^2 = (z-A)^{-1}, \quad z \in \rho(A).
\end{equation}
Then the operator $T$ reads
\begin{equation}\label{T.def}
	T = A^{\frac 12}(I-B(0))A^{\frac 12},
\end{equation}
where $B(z)$, $z \in \rho(A)$, is a bounded operator uniquely determined by the bounded form 
$v(K(z) \cdot,K(z)^* \cdot)$; i.e.~$B(z)$ satisfies
\begin{equation}\label{B.form.def}
	\langle B(z)f,g \rangle = v(K(z) f,K(z)^* g), \qquad f,g \in \cH. 
\end{equation}
In fact, the operator $B$ is the precise meaning of the term $KVK$ in \eqref{res.repr.intro} and elsewhere in the introduction. 

For all $z \in \rho(A)$, we have
$
	z-T = K(z)^{-1}(I-B(z)) K(z)^{-1},
$
which yields the usual factorization of the resolvent of $T$, namely
\begin{equation}\label{Tz.res.dec}
	(z-T)^{-1} = K(z)(I-B(z))^{-1}K(z),
\end{equation}
provided $I-B(z)$ is invertible and $z \in \rho(A)$; see also~\cite[Lemma~1]{Agranovich-1994-28}, \cite[Chap.VI.3.1]{Kato-1966}. Finally, notice that
\begin{equation}\label{B.HS}
	\begin{aligned}
		\|B(z)\|^2 &\leq \|B(z)\|_{\rm HS}^2 = \sum_{j,k=1}^\infty |\langle B(z) \psi_j,\psi_k \rangle|^2 = 
		\sum_{j,k=1}^\infty \frac{|v(\psi_j,\psi_k)|^2}{|z-\mu_j||z-\mu_k|} 
		\\
		& \leq \left(\sum_{j=1}^\infty \frac{\omega_j^2}{|z-\mu_j|}\right)^2.
	\end{aligned}
\end{equation}
%

Under an additional condition on the $p$-summability of $\{1/\mu_n\}_{n \in \N}$, we find that the eigensystem of $T$ is complete.

\begin{proposition}\label{prop:T.complete}
	Let the operator $A$ be as in \eqref{A.def}, let conditions \eqref{asm:v.loc.sub} and \eqref{asm:om.l1} hold and let $T$ be the m-sectorial operator defined in Proposition~\ref{prop:T.basic.def}. Suppose in addition that there exists $p>0$ such that
	\begin{equation}\label{asm:mu.n.Schatten}
		\sum_{n =1}^\infty \frac{1}{\mu_n^p} < \infty;
	\end{equation}
	i.e.~the resolvent of $A$ belongs to the Schatten class $\cS_p$. Then resolvent of $T$ is in the Schatten class $\cS_p$. Moreover,  the system of Riesz projections $Q_{\lambda}$, $\la \in \sigma(T)$, of $T$ is complete.
\end{proposition}

\subsection{Main result (Riesz basis property)}
\label{ssec:main.res}

We denote the gaps of $\mu_n$, $n \in \N$, by
\begin{equation}\label{gn.def}
	r_n^+:=\mu_{n+1}-\mu_{n}, \quad  
	r_n^-:= \mu_n-\mu_{n-1}, \quad n>1,
\end{equation}
and introduce 
\begin{equation}\label{rn.def}
	r_n := \frac 12 \min\{r_n^-, r_n^+\} = \frac12 \dist(\mu_n, \sigma(T) \setminus \{\mu_n\}), \quad n>1.
\end{equation}
For $N,h_1,h_2>0$, let 
\begin{equation}\label{Pi.def}
	\begin{aligned}
		\Pi_0 &\equiv \Pi_0(N,h_1,h_2):= \left\{z \in \C \, : \,  -h_1 <\Re z \leq \mu_N + r_N, \ |\Im z| \leq h_2 \right\}, 
		\\ \Pi_k&:= B_{r_k}(\mu_k), \quad 
		\Pi:=\bigcup_{k \in \N} \Pi_k, \quad \Gamma_k :=\partial \Pi_k, \quad k \in \N.
	\end{aligned}
\end{equation}
Let $S_0^0$ and $P_k^0$, $k \in \N$, be the Riesz projections of the unperturbed operator $A$ corresponding to the eigenvalues lying in $\Pi_0$ and $\Pi_k$, respectively; i.e. 
\begin{equation}\label{Pn.SN.0.def}
	\begin{aligned}
		S_0^0:=\frac{1}{2\pi \ii} \int_{\partial \Pi_0}(z-A)^{-1} \dd z, \quad 
		P_k^0  := \frac{1}{2\pi \ii} \int_{\Gamma_k}(z-A)^{-1} \dd z, \quad k \in \N.
	\end{aligned}
\end{equation}

We already know from Proposition~\ref{prop:T.complete} that under the assumptions \eqref{asm:v.loc.sub}, \eqref{asm:om.l1} and \eqref{asm:mu.n.Schatten}, $T$ has compact resolvent in $\cS_p$ and the system of the Riesz projections $Q_\la$, $\la \in \sigma(T)$, of $T$ is complete. Now we employ the additional (main) condition 
\begin{equation}\label{asm:omega.G}
	\sum_{\substack{j=1 \\j \neq n}}^\infty \frac{\omega_j^2}{|\mu_n-\mu_j|} + \frac{\omega_n^2}{r_n} = o(1), \quad n \to \infty.
\end{equation}
In the first step, we localize large eigenvalues of $T$.  
\begin{proposition}\label{prop:localization}
Let the operator $A$ be as in \eqref{A.def} and let conditions \eqref{asm:v.loc.sub}, \eqref{asm:om.l1}, \eqref{asm:mu.n.Schatten} and  \eqref{asm:omega.G} hold. Let $T$ be the m-sectorial operator defined in Proposition~\ref{prop:T.basic.def}. 
Then there exist $N_0 \in \N$ and $h_1, h_2>0$ such that 
	\begin{equation}\label{T.sp.loc}
		\sigma(T) = \spp(T) \subset \Pi_0(N_0,h_1,h_2) \cup \bigcup_{k > N_0} \Pi_k;
	\end{equation}
	see \eqref{Pi.def}. With these $N_0$ and $h_1$, $h_2$, Riesz projections 
	\begin{equation}\label{Pn.SN.def}
		\begin{aligned}
			S_0:=\frac{1}{2\pi \ii} \int_{\partial \Pi_0}(z-T)^{-1} \dd z, \quad 
			P_n  := \frac{1}{2\pi \ii} \int_{\Gamma_n}(z-T)^{-1} \dd z, \quad n> N_0,
		\end{aligned}
	\end{equation}
	are well-defined and 
	\begin{equation}\label{rank.SN.Pn}
		\begin{aligned}
			\Rank  S_0 & = \Rank \, S_0^0 = N_0, 
			\\
			\Rank  P_n & = \Rank \,  P_n^0 = 1, \quad  n>N_0.
		\end{aligned}
	\end{equation}
	Moreover, the system of projections
	\begin{equation}\label{S0Pk.syst}
	\{S_0\} \cup \{P_k\}_{k>N_0}	
	\end{equation}
is disjoint and complete.
\end{proposition}

Finally, we show the Riesz property of spectral projections of $T$.

\begin{theorem}\label{thm:RB.slow}
Let the operator $A$ be as in \eqref{A.def} and let conditions \eqref{asm:v.loc.sub}, \eqref{asm:om.l1}, \eqref{asm:mu.n.Schatten} and  \eqref{asm:omega.G} hold. Let $T$ be the m-sectorial operator defined in Proposition~\ref{prop:T.basic.def}. 
Then there exists a boundedly invertible $W \in \cB(\cH)$ such that the system of projections \eqref{S0Pk.syst} satisfies
\begin{equation}\label{S0Pk.Riesz}
S_0 = W^{-1} S_0^0 W, \qquad P_k = W^{-1} P_k^0 W, \quad k > N_0.
\end{equation}
\end{theorem}
\begin{remark}
The claim \eqref{S0Pk.Riesz} implies that it is possible to find a Riesz basis $\{\phi_k\}_{k \in \N}$ where $\phi_k \in \Ran(P_k)$, $k > N_0$, are eigenvectors of $T$ and $\phi_k \in \Ran(S_0)$, $k = 1,\dots,N_0$ are eigenvectors or root vectors (generalized eigenvectors) of $T$. 
\end{remark}

\section{Proofs}
\label{sec:proofs}

\subsection{Perturbed operator $T$ and completeness of its eigensystem}
\begin{proof}[Proof of Proposition~\ref{prop:T.basic.def}]
	Let $z_0 \geq 0$ and let $f \in \cD_a$. We expand $f$ in the orthonormal basis $\{\psi_j\}_{j\in \N}$ of $\cH$,
	\begin{equation}\label{f.Da.exp}
		f = \sum_j f_j \psi_j, \quad \|f\|^2 = \sum_j|f_j|^2, \quad \|f\|_{\cD_a}^2 = a[f] = \|A^\frac12 f\|^2 =  \sum_j \mu_j|f_j|^2.
	\end{equation}
	Then, using \eqref{asm:v.loc.sub}, we obtain
	\begin{equation}\label{v.bdd}
		\begin{aligned}
			|v[f]| &\leq \sum_{j,k} |f_j||f_k| |v(\psi_j,\psi_k)| 
			\leq 
			\Big(\sum_{j} \omega_j |f_j|  \Big)^2			
			\\	&=
			\Big(\sum_{j} \frac{\omega_j}{(\mu_j+z_0)^\frac12} (\mu_j+z_0)^\frac12 |f_j|  \Big)^2
			\\
			&\leq 
			\Big(\sum_j \frac{\omega_j^2}{\mu_j+z_0}
			\Big)
			\Big(
			\|f\|_{\cD_a}^2 + z_0 \|f\|^2
			\Big).
		\end{aligned}
	\end{equation}
	Employing the condition \eqref{asm:om.l1} and setting $z_0 = 0$, we find that the quadratic form $v$ is bounded on $\cD_a$; the boundedness of the sesquilinear form follows by a standard polarization argument; see e.g.~\cite[Lemma~IV.2.1]{EE}. Next, by the dominated convergence, the sum on the r.h.s.~of \eqref{v.bdd} decays as $z_0 \to +\infty$, hence \eqref{v.eps.a} holds.  The compactness of the resolvent follows from \eqref{v.eps.a} and \cite[Thm.~VI.3.4]{Kato-1966}.
\end{proof}

\begin{proof}[Proof of Proposition~\ref{prop:T.complete}]
	We verify the assumptions of \cite[Cor.XI.9.31]{DS2}, recalled as Theorem~\ref{thm:compl} above. 
	
	First, since $A^{-1} \in \cS_p$, we have $K(z) \in \cS_{2p}$, $ z \in \rho(A)$. Moreover, $\|B(z_0)\|<1$ for $z_0<0$ with a sufficiently large $|z_0|$ since \eqref{B.HS} and \eqref{asm:om.l1} together with the dominated convergence yield
	\begin{equation}
	\|B(z_0)\| \leq \sum_{j=1}^\infty \frac{\omega_j^2}{|z_0|+\mu_j} \to 0, \quad z_0 \to -\infty.
	\end{equation}
	Thus by the resolvent representation \eqref{Tz.res.dec}, we find that $(z_0-T)^{-1} \in \cS_p$.
	
	Next, we show below that, for every $\vartheta \in [-\pi,\pi]$ with $\vartheta \neq 0$, 
	\begin{equation}\label{T.res.ray}
		\|(e^{\ii \vartheta} \rho - T)^{-1}\| = \BigO(\rho^{-1}), \quad \rho \to + \infty,	
	\end{equation}
	and thus \cite[Cor.XI.9.31]{DS2} indeed yields the claimed completeness.
	
	To establish \eqref{T.res.ray}, let $\vartheta \in (0,\pi]$, $\eps \in (0,1)$ and let $C_\eps>0$ be the constant from \eqref{v.eps.a}. Then for any $f\in \cD_a$
	\begin{equation}
	|\Im t[f]| =  |\Im (a[f] + v[f])| = |\Im v[f]| \leq |v[f]|,
	\end{equation}
	and using \eqref{v.eps.a} in the last step,
	\begin{equation}
	\begin{aligned}
	\Re t[f] & = \Re(a[f]+v[f]) = a[f] + \Re v[f] \geq a[f] - |v[f]|
	\\
		& \geq a[f] - \eps(a[f] + C_\eps \|f\|^2).
	\end{aligned}
	\end{equation}
	Hence, we arrive at
	\begin{equation}
		\frac{|\Im t[f]|}{\Re t[f] + C_\eps\|f\|^2} \leq 
		\frac{|v[f]|}{a[f] - \eps(a[f] + C_\eps\|f\|^2)+ C_\eps\|f\|^2} \leq \frac{\eps}{1-\eps}.
	\end{equation}
	Thus we can enclose the numerical range of $T+ C_\eps$ 
	\begin{equation}
	\Num(T+C_\eps) = \{ \langle (T+C_\eps)f,f\rangle \, : \, f \in \Dom(T), \|f\|=1\},	
	\end{equation}
	see \cite[Chap.~V.3.2, VI.1.2]{Kato-1966}, \cite[Sect.~26.1]{Markus-1988}, in a sector with the semiangle of size $\BigO(\eps)$, $\eps \to 0+$. More precisely,
	\begin{equation}
	\Num(T+C_\eps) \subset \Num (t+C_\eps) 
	 \subset \big\{ z \in \C \, : \, \Re z \geq 0, |\arg z| \leq  \arctan \tfrac{\eps}{1-\eps} \big \}.	
	\end{equation}
	For an arbitrary $\vartheta \in (0,\pi]$ take $\eps = \vartheta / (2 + \vartheta)$. Then
	\begin{equation}\label{eps.vartheta}
	\arctan \frac{\eps}{1-\eps} =  \arctan \frac \vartheta 2 < \frac \vartheta 2, 
	\end{equation}
	from which we obtain
	\begin{equation}
	\Num(T+C_\eps) \subset \Sigma_{\frac \vartheta 2} := \{  z \in \C \, : \, \Re z \geq 0, |\arg z| \leq \tfrac \vartheta 2\}.
	\end{equation}
	Thus
	
	\begin{equation}
	\dist (e^{\ii \vartheta} \rho, \Num(T+C_\eps)) \geq \dist (e^{\ii \vartheta} \rho, \Sigma_{\frac \vartheta2}) = \rho \sin \tfrac \vartheta 2
	\end{equation}
	and, by \cite[Thm.~V.3.2]{Kato-1966} or \cite[Thm.~26.2]{Markus-1988}, 
	\begin{equation}
	\|(e^{\ii \vartheta} \rho -C_\eps - T)^{-1}\| \leq \frac{1}{\dist (e^{\ii \vartheta} \rho, \Num(T+C_\eps))}
	\leq \frac{1}{\rho \sin \frac \vartheta 2}.
	\end{equation}
	Finally, the first resolvent identity yields

	\begin{equation}
		\begin{aligned}
			\|(e^{\ii \vartheta} \rho - T)^{-1}\| &\leq
			\|(e^{\ii \vartheta} \rho -C_\eps - T)^{-1}\|(1+ C_\eps \|(e^{\ii \vartheta} \rho - T)^{-1}\|)
			\\
			& \leq \frac{1}{\rho \sin \frac \vartheta 2 } (1+ C_\eps \|(e^{\ii \vartheta} \rho - T)^{-1}\|),
		\end{aligned}	 
	\end{equation}
and if 
	\begin{equation}
	\rho > \frac{2 C_\eps}{ \sin \frac \vartheta 2}, 
	\end{equation}
	then 
\begin{equation} 
	\|(e^{\ii \vartheta} \rho - T)^{-1}\| \leq \frac 2 {\rho \sin \frac \vartheta2},
\end{equation}
	so \eqref{T.res.ray} follows. 
	With \eqref{eps.vartheta} and $\eps \in (0,1)$ being arbitrary, we verified \eqref{T.res.ray} for any $\vartheta \in (0,\pi]$.
	The argument for $\vartheta \in [-\pi,0)$ is completely analogous.
\end{proof}

\subsection{Localization of the eigenvalues of the perturbed operator}
\label{ssec:loc}
In this section we prove Proposition~\ref{prop:localization}.

Let %
\begin{equation}\label{Sigma.def}
\Xi_n:= \left\{ z \in \C \, : \, \mu_n - \frac 12 r_n^- \leq \Re z \leq \mu_n + \frac 12 r_n^+ \right\}, \quad n \in \N; 
\end{equation}
i.e.~the vertical strips around the eigenvalues $\mu_n$. Further, let
\begin{equation}\label{sigma.N.def}
	\sigma_N:= \sup_{n \geq N} \left(\sum_{\substack{j=1 \\j \neq n}}^\infty \frac{\omega_j^2}{|\mu_n-\mu_j|} + \frac{\omega_n^2}{r_n} \right), \quad N \in \N.
\end{equation}

\begin{lemma}\label{lem:sum.KVK.norm}
Let $\{\omega_j\}_{j \in \N}$ satisfy \eqref{asm:om.l1} and \eqref{asm:omega.G}. Let $\Pi$ be as in \eqref{Pi.def} and $\sigma_N$ as in \eqref{sigma.N.def}. Then 
\begin{equation}\label{sum.KVK.norm.Re}
\sigma_N' :=	\sup_{\substack{z \notin \Pi \\[1mm] \Re z \geq \mu_N - \frac 12 r_N^-}}   \sum_{j=1}^\infty \frac{\omega_j^2}{|z-\mu_j|} \leq 2  \sigma_N, \quad N>1.
\end{equation}	
Moreover, for any $N \in \N$ fixed,
\begin{equation}\label{sum.KVK.norm.Im}
\sup_{\substack{z \notin \Pi \\[1mm] |\Re z| \leq \mu_N, |\Im z| \geq Y} }  \sum_{j=1}^\infty \frac{\omega_j^2}{|z-\mu_j|} = o(1), \quad Y \to +\infty,
\end{equation}
and 
\begin{equation}\label{sum.KVK.norm.neg}
\sup_{\Re z \leq -h  }  \sum_{j=1}^\infty \frac{\omega_j^2}{|z-\mu_j|} = o(1), \quad h \to +\infty. 
\end{equation}
\end{lemma}
\begin{proof}
To justify \eqref{sum.KVK.norm.Re}, let $N \in \N \setminus\{1\}$ be fixed and let $z \notin \Pi$ satisfy $z \in \Xi_n$ for some $n \in \N$, $n \geq N$. First, since $|z-\mu_n| \geq r_n$,
\begin{equation}\label{sum.KVK.z.1}
\sum_{j} \frac{\omega_j^2}{|z-\mu_j|} \leq \sum_{j \neq n} \frac{\omega_j^2}{|\Re z-\mu_j|} + \frac{\omega_n^2}{|z-\mu_n|} \leq 
\sum_{j \neq n} \frac{\omega_j^2}{|\Re z-\mu_j|} + \frac{\omega_n^2}{r_n}.
\end{equation} 
Next, for $j < n$,
\begin{equation}\label{sum.KVK.z.2a}
\begin{aligned}
|\Re z-\mu_j| & = \Re z - \mu_j \geq \mu_n - \frac 12 r_n^- - \mu_j = \mu_n - \mu_j - \frac{\mu_n-\mu_{n-1}}{2}
\\
& = \frac 12  (\mu_n-\mu_j) + \frac 12 (\mu_{n-1} - \mu_j) \geq \frac 12  (\mu_n-\mu_j),
\end{aligned}
\end{equation}
and similarly, for $j > n$,
\begin{equation}\label{sum.KVK.z.2b}
	\begin{aligned}
		|\Re z-\mu_j| & =  \mu_j - \Re z \geq \mu_j - \mu_n - \frac 12 r_n^+ = \mu_j - \mu_n - \frac{\mu_{n+1}-\mu_{n}}{2}
		\\
		& = \frac 12  (\mu_j-\mu_n) + \frac 12 (\mu_j - \mu_{n+1}) \geq \frac 12  (\mu_j-\mu_n).
	\end{aligned}
\end{equation}
Hence, for any $n \geq N$ and any $z \in \Xi_n$ with $z \notin \Pi$, 
\begin{equation}\label{sum.KVK.z.3}
\sum_{j} \frac{\omega_j^2}{|z-\mu_j|} \leq 2 \sum_{j \neq n} \frac{\omega_j^2}{|\mu_n-\mu_j|} + \frac{\omega_n^2}{r_n} \leq 2 \sigma_N;
\end{equation} 
i.e.~\eqref{sum.KVK.norm.Re} holds.

To show \eqref{sum.KVK.norm.Im}, let $z \notin \Pi$ and $|\Re z| \leq \mu_N$ with $N \in \N$ arbitrary but fixed and let $|\Im z | > Y$. Let $\tilde N \in \N$ be such that
\begin{equation}\label{N.tild.def}
\mu_{\tilde N} \geq 2 \mu_N
\end{equation}
and we split the summation as
\begin{equation}\label{sum.split}
\sum_{j} \frac{\omega_j^2}{|z-\mu_j|} = \sum_{j \leq \tilde N } \frac{\omega_j^2}{|z-\mu_j|} + \sum_{j > \tilde N } \frac{\omega_j^2}{|z-\mu_j|}.
\end{equation}
First, the finite sum can be estimated straightforwardly as
\begin{equation}\label{fin.sum.est}
\sum_{j \leq \tilde N} \frac{\omega_j^2}{|z-\mu_j|}  
\leq 
\frac{1}{Y} \sum_{j \leq \tilde N} \omega_j^2.
\end{equation}
To estimate the second sum, notice that for $j > \tilde N$
\begin{equation}
|\Re z - \mu_j| \geq \mu_j - \mu_N \geq \frac{\mu_j}2  +  \frac{\mu_{\tilde N}}2  - \mu_N \geq \frac{\mu_j}2; 
\end{equation}
see \eqref{N.tild.def}.
Thus
\begin{equation}
\sum_{j > \tilde N } \frac{\omega_j^2}{|z-\mu_j|}
\leq 2 \sum_{j > \tilde N } \frac{\omega_j^2}{|\Re z-\mu_j| + |\Im z|} 
\leq 4 \sum_{j > \tilde N } \frac{\omega_j^2}{\mu_j + Y}.
\end{equation}
Since the sequence $\{\omega_j^2/\mu_j\}_{j \in \N}$ is summable by the assumption \eqref{asm:om.l1}, the dominated convergence theorem yields that
\begin{equation}\label{inf.sum.est}
\sum_{j > \tilde N } \frac{\omega_j^2}{|z-\mu_j|} = o(1), \quad Y \to + \infty,
\end{equation}
hence, together with \eqref{fin.sum.est}, we verified that \eqref{sum.KVK.norm.Im} holds.

Finally, to establish \eqref{sum.KVK.norm.neg}, let $h>0$ and notice that for $z$ with $\Re z \leq -h$, 
\begin{equation}
\sum_{j} \frac{\omega_j^2}{|z-\mu_j|} \leq  \sum_{j} \frac{\omega_j^2}{h+\mu_j}.
\end{equation}
Thus \eqref{sum.KVK.norm.neg} follows again by the dominated convergence.
\end{proof}

\begin{proof}[Proof of Proposition~\ref{prop:localization}]
The key step of the proof is the estimate of the norm of $B(z)$; see Section~\ref{ssec:T.def} and in particular recall \eqref{B.HS}, namely
	\begin{equation}\label{KVK.HS}
		\|B(z)\| \leq \sum_{j} \frac{\omega_j^2 }{|z-\mu_j|}.
	\end{equation}
	First, from \eqref{sum.KVK.norm.neg} in Lemma~\ref{lem:sum.KVK.norm} we obtain 
	\begin{equation}
		\sup_{\Re z \leq -h} \|B(z)\| = o(1), \quad h \to +\infty. 
	\end{equation}
	Thus there exists $h_1>0$ such that
	\begin{equation}
		\sup_{\Re z \leq -h_1} \|B(z)\| \leq \frac 12. 
	\end{equation}
	
	Next, from \eqref{sum.KVK.norm.Re}, 
	\begin{equation}
	\sup_{\substack{z \notin \Pi \\[1mm] \Re z \geq \mu_N-r_N}}\|B(z)\|	\leq \sup_{\substack{z \notin \Pi \\[1mm] \Re z \geq \mu_N-r_N^-}}\|B(z)\| = \BigO(\sigma_N), \quad N \to + \infty.
	\end{equation}
	Thus, from the assumption \eqref{asm:omega.G}, there exists $N_0 \in \N$ such that
	\begin{equation}\label{KVK.1/2}
		\sup_{\substack{z \notin \Pi \\[1mm] \Re z \geq \mu_{N_0} - r_{N_0}}}\|B(z)\| \leq \frac12.
	\end{equation}
	Finally, with these $N_0$ and $h_1$ fixed, \eqref{sum.KVK.norm.Im} yields
	\begin{equation}
		\sup_{\substack{-h_1 \leq \Re z \leq \mu_{N_0} + r_{N_0} \\[1mm] |\Im z| \geq h }} \|B(z)\| \leq 
		\sup_{\substack{-h_1 \leq \Re z \leq \mu_{N_0+1} \\[1mm] |\Im z| \geq h }} \|B(z)\| = o(1), \quad h \to + \infty;
	\end{equation}
	hence there exists $h_2>0$ such that
	\begin{equation}
		\sup_{\substack{-h_1 \leq \Re z \leq \mu_{N_0} + r_{N_0} \\[1mm] |\Im z| \geq h_2 }} \|B(z)\| \leq \frac 12.
	\end{equation}
	The norm estimates of $B(z)$ above and the resolvent representation \eqref{Tz.res.dec} justify the localization of the eigenvalues in \eqref{T.sp.loc}. Moreover, the spectral projections in \eqref{Pn.SN.def} are well-defined and disjoint; see Lemma~\ref{lem:disj}.
	
	The claim on the rank \eqref{rank.SN.Pn} follows by a standard argument based on that 
	\begin{equation}
		t \mapsto \Tr \, \frac{1}{2\pi \ii} \int_{\Gamma_n} K(z)(I-t B(z))^{-1} K(z) \, \dd z, \quad 0 \leq t \leq 1,   
	\end{equation}
	is a continuous integer-valued function, thus this function is constant; see also \cite[Lemma~VII.6.7]{DS1}.

	Finally, the completeness follows from Proposition~\ref{prop:T.complete}.
\end{proof}

\subsection{Proof of Theorem~\ref{thm:RB.slow}}
\label{ssec:proof.thm}

One of the main steps in the proof is the norm estimate of the operators $\cM$, $\cM'$ defined as follows. Let $N \in \N$,  
\begin{equation}\label{chiN.def}
	\chi_N(n):= \begin{cases}
		0, & 0<n < N,
		\\
		1 & n \geq N,
	\end{cases} 
\end{equation}
and let $\cM$ be the operator in $\ell^2(\N)$ with matrix elements
\begin{equation}\label{cM.def}
	\cM_{nk} : = 
	\begin{cases}
		\displaystyle
		\chi_{N}(n) \frac{\omega_n \omega_k}{|\mu_n-\mu_k|}, & n \neq k, 
		\\[2mm]
		0, & n=k, \qquad k,n \in \N.	
	\end{cases}
\end{equation}
Moreover, let $z_n \in \Gamma_n$, $n \in \N$, and let $\cM'$ be the operator in $\ell^2(\N)$ with matrix elements
\begin{equation}\label{cM'.def}
\cM_{nk}' : = \chi_{N}(n) \frac{\omega_n \omega_k }{|z_n-\mu_k|}, \qquad k,n \in \N.
\end{equation}
Notice that $\cM \equiv \cM(N)$, i.e.,~it depends on $N$, and $\cM' \equiv \cM'(N,\{z_n\}_{n \in \N})$, i.e.,~it depends on $N$ and the sequence $\{z_n\}_{n \in \N}$.
	
\begin{lemma}\label{lem:schur}
Let the assumptions of Theorem~\ref{thm:RB.slow} be satisfied, let $N \in \N$, let $\sigma_N$ be as in \eqref{sigma.N.def}, let
\begin{equation}\label{rho.def}
\varrho_N := \sum_{n \geq N} \frac{\omega_n^2}{\mu_n},  \quad N \in \N,
\end{equation}
and let $k_N \in \N$ be the largest $k_N \leq N$ such that 
\begin{equation}\label{kN.def}
\mu_N -r_N \geq 2 \mu_{k_N};
\end{equation}
such $k_N$ exists for all $N > \wt N$ with some $\wt N \in \N$ sufficiently large, moreover $k_N \to \infty$ and $N \to \infty$.
Further, let $\cM$, $\cM'$ be the operators as in \eqref{cM.def}, \eqref{cM'.def}, respectively. 
Then, for all $N \in \N$ with $N > \wt N$,
\begin{equation}\label{cM.est}
\|\cM\|^2 \leq \sigma_N \max\{2\varrho_N, \sigma_{k_N}\}, \quad 
\|\cM'\|^2 \leq 4 \sigma_N \max\{2\varrho_N,  \sigma_{k_N}\}.
\end{equation}

\end{lemma}
\begin{proof}
Both estimates in \eqref{cM.est} follow by the Schur test with the weight $\{\omega_k\}_{k \in \N}$; see e.g.~\cite{Schur-1911-140}, \cite[Thm.5.2]{Halmos-1978}, \cite[Sec.~3]{Dym-2003-210}.

First for $\cM$, we obtain
\begin{align}
\sum_k \cM_{nk} \omega_k & \leq \chi_N(n) \omega_n \sum_{k \neq n} \frac{\omega_k^2}{|\mu_n-\mu_k|} \leq \sigma_N \omega_n, \quad n \in \N,
\label{Schur.1}
\\
\sum_n \cM_{nk} \omega_n & \leq \omega_k \sum_{\substack{n \neq k \\ n \geq N} } \frac{\omega_n^2}{|\mu_n-\mu_k|}, \quad k \in \N. \label{Schur.2}
\end{align}	
To estimate the sum on the r.h.s.~of \eqref{Schur.2}, we consider two cases. If $k < k_N$, then for $n \geq N$, 
\begin{equation}
\mu_n-\mu_k \geq  \frac{\mu_n}2 + \frac{\mu_N}2 - \mu_{k_N} \geq \frac{\mu_n}2; 
\end{equation}
see~\eqref{kN.def}. 
Hence, in this case ($k < k_N$),
\begin{equation}\label{cM.small.k}
\sum_{\substack{n \neq k \\ n \geq N} } \frac{\omega_n^2}{|\mu_n-\mu_k|}  
\leq 
2 \sum_{n \geq N} \frac{\omega_n^2}{\mu_n} = 2 \varrho_N.
\end{equation}
On the other hand, if $k \geq k_N$, then
\begin{equation}\label{cM.large.k}
\sum_{\substack{n \neq k \\ n \geq N} } \frac{\omega_n^2}{|\mu_n-\mu_k|}
\leq
\sum_{n \neq k} \frac{\omega_n^2}{|\mu_n-\mu_k|}
\leq 
\sigma_{k_N}.
\end{equation}
Thus, from \eqref{Schur.2}, \eqref{cM.small.k} and \eqref{cM.large.k}, we arrive at
\begin{equation}\label{Schur.2.fin}
\sum_n \cM_{nk} \omega_n \leq \max\{2\varrho_N,\sigma_{k_N}\} \omega_k
\end{equation}
and, employing \eqref{Schur.1}, the claim on $\|\cM\|$ in \eqref{cM.est} follows by the Schur test.

The case of $\|\cM'\|$ is analogous, in addition, we use \eqref{sum.KVK.norm.Re} and arguments in its proof; see \eqref{sum.KVK.z.1} -- \eqref{sum.KVK.z.3}. In detail, 
\begin{align}
	\sum_k \cM_{nk}' \omega_k & \leq \chi_N(n) \omega_n \sum_{k} \frac{\omega_k^2}{|z_n-\mu_k|} 
	\leq 
	\omega_n \sup_{\substack{z \notin \Pi \\ \Re z \geq \mu_N - r_N}}\sum_{k} \frac{\omega_k^2}{|z-\mu_k|}
	\\
	&
	\leq 2 \sigma_N \omega_n, \quad n \in \N,
	\label{Schur.3}
	\\
	\sum_n \cM_{nk}' \omega_n & \leq \omega_k \sum_{n \geq N } \frac{\omega_n^2}{|z_n-\mu_k|}, \quad k \in \N. \label{Schur.4}
\end{align}	
Similarly as above, let first $k <k_N$, thus $\mu_k < \mu_{k_N} \leq (\mu_N-r_N)/2$. Then, for $n \geq N$ and $z_n \in \Gamma_n$, we have $\Re z_n \geq \mu_N  - r_N$ and
\begin{equation}
\Re z_n- \mu_k \geq \frac{\Re z_n}2 + \frac{\mu_N-r_N}2   - \mu_{k_N} \geq \frac{\Re z_n}{2} \geq \frac{\mu_n}2 - \frac{\mu_n-\mu_{n-1}}{4} \geq \frac{\mu_n}{4};
\end{equation}
see \eqref{kN.def} and \eqref{rn.def}. Hence, for $k < k_N$,
\begin{equation}\label{cM'.small.k}
	\sum_{n \geq N} \frac{\omega_n^2}{|z_n-\mu_k|}
	\leq
	\sum_{n \geq N} \frac{\omega_n^2}{\Re z_n-\mu_k} 
	\leq
	4 \sum_{n \geq N} \frac{\omega_n^2}{\mu_n} = 4 \varrho_N.
\end{equation}
Let further $k \geq k_N$. Following and adjusting \eqref{sum.KVK.z.1} -- \eqref{sum.KVK.z.3}, we find that
\begin{equation}\label{cM'.large.k}
\sum_{n \geq N} \frac{\omega_n^2}{|z_n-\mu_k|}
\leq
2 \sum_{n \neq k} \frac{\omega_n^2}{|\mu_n-\mu_k|} +  \frac{\omega_k^2}{r_k}
\leq
2 \sigma_{k_N}.
\end{equation}
Combining \eqref{Schur.4}, \eqref{cM'.small.k} and \eqref{cM'.large.k}, we arrive at
\begin{equation}
\sum_n \cM_{nk}' \omega_n \leq 2 \max\{2 \varrho_N, \sigma_{k_N}\}  \omega_k , \quad k \in \N,
\end{equation}
thus, together with \eqref{Schur.3}, the Schur test yields the claim on $\|\cM'\|$ in \eqref{cM.est}.
\end{proof}

\begin{lemma}\label{lem:KVK.s.exp}
Let the assumptions of Theorem~\ref{thm:RB.slow} be satisfied, let $f \in \cH$ and let
\begin{equation}
	f_j:= \langle f, \psi_j \rangle, \qquad  v_{jk}:= v(\psi_j,\psi_k), \quad \forall j,k \in \N.
\end{equation}	
Then
\begin{equation}\label{K2Vk2.exp}
\langle KBK f,f \rangle =  \sum_{j,k} \frac{v_{jk}}{(z-\mu_j)(z-\mu_k)} f_j \ov{f_k}
\end{equation}
and, for $s \in \N$,
\begin{equation}\label{KVK.s.exp}
\begin{aligned}
&\langle KB^{s+1} K f,f \rangle 
\\
& \qquad = \sum_{j,j_1, \dots, j_s, k} \frac{v_{jj_1} v_{j_1 j_2} \dots v_{j_{s-1} j_s} v_{j_s k} }{(z-\mu_j)(z-\mu_{j_1})(z-\mu_{j_2})\dots (z-\mu_{j_s})(z-\mu_k)} f_j \ov{f_k}.  
\end{aligned}
\end{equation}
Moreover, the series on the r.h.s.~of~\eqref{K2Vk2.exp} and \eqref{KVK.s.exp} converge absolutely.
\end{lemma}
\begin{proof}
The formulas follow by expanding $f$ in the orthonormal basis $\{\psi_j\}_{j \in \N}$ and a repeated use of (see \eqref{asm:v.loc.sub} and \eqref{B.form.def})
\begin{equation}
B \psi_l =
\sum_m \langle B \psi_l, \psi_m \rangle \psi_m = 
 \sum_m \frac{v_{lm}}{(z-\mu_l)^\frac12 (z-\mu_m)^\frac12} \psi_m, \quad l \in \N.
\end{equation}

For the absolute convergence, notice that
\begin{equation}\label{KVK.sums.conv}
\begin{aligned}
&\sum_{j,j_1, \dots, j_s, k} \left|\frac{v_{jj_1} v_{j_1 j_2} \dots v_{j_{s-1} j_s} v_{j_s k} }{(z-\mu_j)(z-\mu_{j_1})(z-\mu_{j_2})\dots (z-\mu_{j_s})(z-\mu_k)} f_j \ov{f_k}\right| 
\\
&\quad  \leq
\left(
\sum_m \frac{\omega_m^2}{|z-\mu_m|}
\right)^s 
\left(
\sum_j \frac{\omega_j}{|z-\mu_j|} |f_j|
\right)^2
\\
& \quad
\leq 
\left(
\sum_m \frac{\omega_m^2}{|z-\mu_m|}
\right)^{s+1} 
\sum_j \frac{|f_j|^2}{|z-\mu_j|} 
\leq 
\left(
\sum_m \frac{\omega_m^2}{|z-\mu_m|}
\right)^{s+1} 
\|K(z)f\|^2
\end{aligned}
\end{equation}
and employ the assumption~\eqref{asm:om.l1}; the absolute convergence in \eqref{K2Vk2.exp} follows by taking $s=0$ in \eqref{KVK.sums.conv}, i.e., the summation  over $j, k$ only.
\end{proof}

In the following, we will use the sequence
\begin{equation}\label{tau.N.def}
\tau_N:= \max\{ \|\cM\|, \|\cM'\|, \sigma_N, \sigma_N' \}. 
\end{equation}
Recall that Lemmas~\ref{lem:sum.KVK.norm} and \ref{lem:schur} show that $\tau_N$ decays, more precisely, 
\begin{equation}
\tau_N =o(\sigma_N^\frac12), \quad N \to + \infty. 
\end{equation}

\begin{lemma}\label{lem:s0.est}
Let the assumptions of Theorem~\ref{thm:RB.slow} be satisfied and let $f \in \cH$. Then
\begin{equation}
\sum_{n>N} 
\left|
\frac{1}{2\pi \ii} \int_{\Gamma_n} \langle KBK f,f \rangle \, \dd z
\right|
\leq 2 \tau_N \|f\|^2
\end{equation}

\end{lemma}
\begin{proof}
The Cauchy formula and \eqref{K2Vk2.exp} lead to
\begin{equation}
\frac{1}{2\pi \ii} \int_{\Gamma_n}	\langle KB(z)K f,f \rangle \ \dd z = 
\sum_{k \neq n} \frac{v_{nk}}{\mu_n-\mu_k} f_n \ov{f_k} 
+
\sum_{j \neq n} \frac{v_{jn}}{\mu_n-\mu_j} f_j \ov{f_n}.
\end{equation}
Thus Cauchy-Schwartz inequality, the assumption~\eqref{asm:v.loc.sub} and Lemma~\ref{lem:schur} yield
\begin{align}
& \sum_{n>N} 
\left|
\frac{1}{2\pi \ii} \int_{\Gamma_n} \langle KBK f,f \rangle \, \dd z
\right|
\leq 2 \sum_{n>N} |f_n| \left(\sum_{k \neq n} \frac{\omega_n \omega_k}{|\mu_n-\mu_k|} |f_k| \right) 
\\ 
& \quad\leq 
2 \left(
\sum_{n>N} \left(
\sum_{k \neq n}\frac{\omega_n \omega_k}{|\mu_n-\mu_k|} |f_k|
\right)^2
\right)^\frac 12 \|f\|
\leq 2 \|\cM\| \|f\|^2. 
\qedhere
\end{align}
\end{proof}

To estimate the terms with $s \in \N$, we split the sum over $s+2$ indices in \eqref{KVK.s.exp} to sums $\Sigma^t$ over $s+2-t$ indices with the remaining exactly $t$ indices equal $n$, schematically, 
\begin{equation}\label{Sigma.t.def}
\sum_{j,j_1, \dots, j_s, k} = 
{\sum}^0 
+
{\sum}^1 
+ \dots + {\sum}^{s+2}. 
\end{equation}

We note first that due to the structure of the summand in \eqref{KVK.s.exp}, the contour integral over $\Gamma_n$ is zero for the sums $\Sigma^0$ and $\Sigma^{s+2}$. 
\begin{lemma}\label{lem:S0}
Let the assumptions of Theorem~\ref{thm:RB.slow} be satisfied, let $s \in \N$ and let $f \in \cH$. Then
\begin{equation}
\begin{aligned}
{\sum}^0 \frac{1}{2\pi \ii} \int_{\Gamma_n}  
\frac{v_{jj_1} v_{j_1 j_2} \dots v_{j_s k} }{(z-\mu_j)(z-\mu_{j_1})\dots (z-\mu_{j_s})(z-\mu_k)} f_j \ov{f_k}  
\, \dd z &= 0,
\\
{\sum}^{s+2} \frac{1}{2\pi \ii} \int_{\Gamma_n}  
\frac{v_{jj_1} v_{j_1 j_2} \dots v_{j_s k} }{(z-\mu_j)(z-\mu_{j_1})\dots (z-\mu_{j_s})(z-\mu_k)} f_j \ov{f_k}  
\, \dd z &= 0.
\end{aligned}
\end{equation}
\end{lemma}
\begin{proof}
From the definition of $\Sigma^0$, we have
\begin{equation}
\begin{aligned}
&{\sum}^0 \frac{1}{2\pi \ii} \int_{\Gamma_n}  
\frac{v_{jj_1} v_{j_1 j_2} \dots v_{j_s k} }{(z-\mu_j)(z-\mu_{j_1})\dots (z-\mu_{j_s})(z-\mu_k)} f_j \ov{f_k}  
\, \dd z
\\
&\quad = 
\sum_{j,j_1,\dots,j_s,k \neq n} \frac{1}{2\pi \ii} \int_{\Gamma_n}  
\frac{v_{jj_1} v_{j_1 j_2} \dots v_{j_s k} }{(z-\mu_j)(z-\mu_{j_1})\dots (z-\mu_{j_s})(z-\mu_k)} f_j \ov{f_k}  
\, \dd z.
\end{aligned}
\end{equation}
Since the integrand is analytic in the disc $\Pi_n$, all contour integrals over $\Gamma_n = \partial \Pi_n$ are zero.

Next, the definition of $\Sigma^{s+2}$ and the residue calculus lead to
\begin{equation}
\begin{aligned}
	&{\sum}^{s+2} \frac{1}{2\pi \ii} \int_{\Gamma_n}  
		\frac{v_{jj_1} v_{j_1 j_2} \dots v_{j_s k} }{(z-\mu_j)(z-\mu_{j_1})\dots (z-\mu_{j_s})(z-\mu_k)} f_j \ov{f_k}  
		\, \dd z
		\\
		& \quad = 
		\frac{1}{2\pi \ii} \int_{\Gamma_n}  
		\frac{v_{nn}^{s+1}}{(z-\mu_n)^{s+2}} |f_n|^2   
		\, \dd z = 0.
\end{aligned} \qedhere
\end{equation}
\end{proof}

In the following, we perform the contour integration for the sum $\Sigma^1$ and estimate the resulting sums. 

\begin{lemma}\label{lem:t1.est}
Let the assumptions of Theorem~\ref{thm:RB.slow} be satisfied, let $s \in \N$ and let $f \in \cH$. Then
\begin{equation}\label{S1.est}
\begin{aligned}
&\sum_{n>N} 
\left|
{\sum}^1 \frac{1}{2\pi \ii} \int_{\Gamma_n}  
\frac{v_{jj_1} v_{j_1 j_2} \dots v_{j_s k} }{(z-\mu_j)(z-\mu_{j_1})\dots (z-\mu_{j_s})(z-\mu_k)} f_j \ov{f_k}  
\, \dd z
\right|
\\
&
\quad \leq \binom{s+2}{1} \tau_N^{s+1} \|f\|^2.
\end{aligned}
\end{equation}
\end{lemma}
\begin{proof}
The contour integration results in $s+2$ sums over $s+1$ indices (none of which is equal $n$), namely
\begin{equation}\label{s1.sums}
	\begin{aligned}
		& \sum_{j_1, j_2, \dots, j_s, k \neq n} \frac{v_{nj_1} v_{j_1 j_2} \dots v_{j_s k} }{(\mu_n-\mu_{j_1})(\mu_n-\mu_{j_2})\dots (\mu_n-\mu_{j_s})(\mu_n-\mu_k)} f_n \ov{f_k}  
		\\
		& + \sum_{j, j_2, \dots, j_s, k \neq n} \frac{v_{jn} v_{n j_2} \dots v_{j_s k} }{(\mu_n-\mu_j)(\mu_n-\mu_{j_2})\dots (\mu_n-\mu_{j_s})(\mu_n-\mu_k)} f_j \ov{f_k}  + \dots +
		\\
		& + \sum_{j, j_1, \dots, j_{s-1}, k \neq n} \frac{v_{jj_1} \dots v_{j_{s-1} n}  v_{n k} }{(\mu_n-\mu_j)(\mu_n-\mu_{j_1})\dots (\mu_n-\mu_{j_{s-1}})(\mu_n-\mu_k)} f_j \ov{f_k}  
		\\
		& + \sum_{j,j_1, j_2, \dots, j_s  \neq n} \frac{v_{jj_1} v_{j_1 j_2} \dots v_{j_s n} }{(\mu_n-\mu_j)(\mu_n-\mu_{j_1})\dots (\mu_n-\mu_{j_s})} f_j \ov{f_n}.
	\end{aligned}
\end{equation}
We start with the estimate of the first term in \eqref{s1.sums}. The assumption~\eqref{asm:v.loc.sub}, Cauchy-Schwartz inequality and Lemma~\ref{lem:schur} lead to (see also \eqref{sigma.N.def})
\begin{equation}
\begin{aligned}
&\sum_{n>N}  \left|
\sum_{j_1, j_2, \dots, j_s, k \neq n} \frac{v_{nj_1} v_{j_1 j_2} \dots v_{j_s k} }{(\mu_n-\mu_{j_1})(\mu_n-\mu_{j_2})\dots (\mu_n-\mu_{j_s})(\mu_n-\mu_k)} f_n \ov{f_k}
\right|
\\
&
\leq \sigma_N^s  \sum_{n>N} \sum_{k \neq n} \frac{\omega_n \omega_k}{|\mu_n-\mu_k|} |f_n||f_k|
\leq \sigma_N^s \Big(\sum_{n>N} \Big(\sum_{k \neq n} \frac{\omega_n \omega_k}{|\mu_n-\mu_k|}|f_k|\Big)^2\Big)^\frac12 \|f\|
\\
& 
\leq \sigma_N^s \|\cM\| \|f\|^2 \leq \tau_N^{s+1} \|f\|^2.
\end{aligned}
\end{equation}
The estimate of the last term in \eqref{s1.sums} is completely analogous.

The remaining terms in \eqref{s1.sums} can be all estimated in the same way. In detail for the term with no summation in $j_1$, employing the assumption~\eqref{asm:v.loc.sub} and Lemma~\ref{lem:schur}, we arrive at
\begin{equation}
\begin{aligned}
&\sum_{n>N}  \left|
\sum_{j, j_2, \dots, j_s, k \neq n} \frac{v_{jn} v_{n j_2} \dots v_{j_s k} }{(\mu_n-\mu_j)(\mu_n-\mu_{j_2})\dots (\mu_n-\mu_{j_s})(\mu_n-\mu_k)} f_j \ov{f_k}
\right|
\\
& \leq 
\sigma_N^{s-1}\sum_{n>N} \sum_{j,k \neq n} \frac{\omega_n^2 \omega_j \omega_k }{|\mu_n-\mu_j||\mu_n-\mu_k|}|f_j| |f_k|
\\
&=
\sigma_N^{s-1}\sum_{n>N} \Big(
\sum_{k\neq n} \frac{\omega_n \omega_k}{|\mu_n-\mu_k|}|f_k| 
\Big)^2
\leq \sigma_N^{s-1} \|\cM\|^2 \|f\|^2 \leq \tau_N^{s+1} \|f\|^2.
\end{aligned}
\end{equation}
In total, there are $s+2$ sums in \eqref{s1.sums}, thus the claim \eqref{S1.est} is justified.
\end{proof}

For $\Sigma^t$ with $1<t<s+2$, we estimate the integrals directly.
\begin{lemma}\label{lem:t.est}
Let the assumptions of Theorem~\ref{thm:RB.slow} be satisfied,  let $f \in \cH$, let $s \in \N$ and $1<t<s+2$. Then
\begin{equation}\label{St.est}
	\begin{aligned}
		&\sum_{n>N} 
		\left|
		{\sum}^t \frac{1}{2\pi \ii} \int_{\Gamma_n}  
		\frac{v_{jj_1} v_{j_1 j_2} \dots v_{j_s k} }{(z-\mu_j)(z-\mu_{j_1})\dots (z-\mu_{j_s})(z-\mu_k)} f_j \ov{f_k}  
		\, \dd z
		\right|
		\\
		&
		\quad \leq \binom{s+2}{t} \tau_N^{s+1} \|f\|^2.
	\end{aligned}
\end{equation}
\end{lemma}
\begin{proof}
First we apply triangle inequality
\begin{equation}
\begin{aligned}
&\sum_{n>N} 
\left|
{\sum}^t \frac{1}{2\pi \ii} \int_{\Gamma_n}  
\frac{v_{jj_1} v_{j_1 j_2} \dots v_{j_s k} }{(z-\mu_j)(z-\mu_{j_1})\dots (z-\mu_{j_s})(z-\mu_k)} f_j \ov{f_k}  
\, \dd z
\right|
\\
&
\leq \sum_{n>N} \frac{1}{2 \pi } {\sum}^t \int_{\Gamma_n} \frac{\omega_j \omega_{j_1}^2 \dots \omega_{j_s}^2 \omega_{k}}{|z-\mu_j||z-\mu_{j_1}|\dots|z-\mu_{j_s}||z-\mu_k|} |f_j||f_k|\, |\dd z|.
\end{aligned}
\end{equation}
Each $\Sigma^t$ can be written as a sum of $\binom{s+2}{t}$ sums in $s+2-t$ indices (with the remaining $t$ indices equal $n$). These resulting sums are of three types, namely 

$S_1$: both $j$ and $k$ are equal $n$; 

$S_2$: either $j$ or $k$ is equal $n$;

$S_3$: neither $j$ nor $k$ is equal $n$.

The sum of type $S_1$ does not exceed (see \eqref{sum.KVK.norm.Re} and \eqref{tau.N.def})
\begin{equation}\label{sum.t.i.1}
\begin{aligned}
&\sum_{n>N}  \frac{\omega_n^{2+2(t-2)}}{r_n^t} \frac{1}{2 \pi }  \int_{\Gamma_n} \left(\sum_{l} \frac{\omega_l^2}{|z-\mu_l|}  \right)^{s-(t-2)} \, |f_n|^2 \, |\dd z|
\\
& \quad \leq 
(\sigma'_N)^{s-(t-2)} \sum_{n>N} \frac{\omega_n^{2(t-1)}}{r_n^t}     \frac{1}{2 \pi }  \int_{\Gamma_n}  |f_n|^2 \, |\dd z|
\\
& \quad \leq
(\sigma'_N)^{s-(t-2)}  \sum_{n>N} \left(\frac{\omega_n^2}{r_n}\right)^{t-1} |f_n|^2
\leq (\sigma'_N)^{s-(t-2)} \sigma_N^{t-1} \|f\|^2 
\leq \tau_N^{s+1} \|f\|^2.
\end{aligned}
\end{equation}

Employing similar steps as above, we find that the sum of the type $S_2$ (with $j=n$, the case $k=n$ is fully analogous) does not exceed
\begin{equation}
\begin{aligned}
&\sum_{n>N}  \frac{\omega_n^{2(t-1)}}{r_n^t}  \frac{1}{2 \pi }  \int_{\Gamma_n} \left(\sum_{l} \frac{\omega_l^2}{|z-\mu_l|}  \right)^{s-(t-1)} \sum_{k} \frac{\omega_n \omega_k}{|z-\mu_k|} |f_n| |f_k|  \, |\dd z| 
\\
& \quad \leq 
(\sigma'_N)^{s-(t-1)}
\sum_{n>N} \frac{\omega_n^{2(t-1)}}{r_n^t}   \frac{1}{2 \pi }  \int_{\Gamma_n} \sum_{k} \frac{\omega_n \omega_k}{|z-\mu_k|} |f_n| |f_k|  \, |\dd z|.
\end{aligned}
\end{equation}
Selecting the points $z_n$, $n >N$, such that 
\begin{equation}\label{zn.sel}
\max_{z \in \Gamma_n}\sum_{k} \frac{\omega_k}{|z-\mu_k|} |f_k| 
= 
\sum_{k} \frac{\omega_k}{|z_n-\mu_k|} |f_k|,
\end{equation}
we infer that the sum of type $S_2$ does not exceed (see Lemma~\ref{lem:schur} and \eqref{tau.N.def})
\begin{equation}
\begin{aligned}
& 
(\sigma'_N)^{s-(t-1)}
\sum_{n>N} \left(\frac{\omega_n^{2}}{r_n}\right)^{t-1} 
\sum_{k} \frac{\omega_n \omega_k}{|z_n-k|} |f_k| |f_n|
\leq \tau_N^{s} \|\cM'\| \|f\|^2 
\\
&
\quad \leq \tau_N^{s+1} \|f\|^2.
\end{aligned}
\end{equation}

Finally, the sum of the type $S_3$ does not exceed
\begin{equation}
\begin{aligned}
&\sum_{n>N}  \frac{\omega_n^{2(t-1)}}{r_n^t} \frac{1}{2 \pi }  \int_{\Gamma_n} \left(\sum_{l} \frac{\omega_l^2}{|z-\mu_l|}  \right)^{s-t} \sum_{j,k} \frac{\omega_n^2 \omega_j \omega_k }{|z-\mu_j||z-\mu_k|} |f_j||f_k|  \, |\dd z| 
\\
& \leq 
\tau_N^{s-t}  \sum_{n>N}   \left(\frac{\omega_n^{2}}{r_n}\right)^{t-1}  
\left(\sum_k \frac{\omega_n \omega_k }{|z_n-\mu_k|} |f_k| \right)^2
\leq 
\tau_N^{s-1}  \|\cM'\|^2 \|f\|^2
\\
& 
\leq 
\tau_N^{s+1} \|f\|^2;
\end{aligned}
\end{equation}
in the first inequality we selected points $\{z_n\}_{n \in \N}$ similarly as in \eqref{zn.sel}.

Since we obtained the same upper bound for the sum of each type and there are $\binom{s+2}{t}$ sums in total, we arrive at the claim \eqref{St.est}.
\end{proof}

\begin{lemma}\label{lem:KVK.s.est}
Let the assumptions of Theorem~\ref{thm:RB.slow} be satisfied, let $f \in \cH$ and let $N_0 \in \N$ be as in Proposition~\ref{prop:localization}. Then there exists $N_* \in \N$ with $N_* > N_0$ such that for each $s \in \N_0$
\begin{equation}\label{N*.s.sum}
\sum_{n \geq N_*} \left| \frac{1}{2 \pi \ii} \int_{\Gamma_n} \langle KB^{s+1} K f,f \rangle \; \dd z \right|
\leq 2^{-s}\|f\|^2.
\end{equation}
\end{lemma}
\begin{proof}
Recall that $\tau_N=o(1)$ as $N \to \infty$, hence the claim for $s=0$ follows immediately from Lemma~\ref{lem:s0.est}; it suffices to choose $N_* > N_0$ such that 
\begin{equation}\label{tau.s0.cond}
\tau_{N_*} \leq \frac 12.
\end{equation}

For $s>0$, we employ Lemma~\ref{lem:KVK.s.exp}, the split of summation as in \eqref{Sigma.t.def} and Lemmas \ref{lem:S0}, \ref{lem:t1.est} and \ref{lem:t.est}. 
In detail
\begin{equation}\label{KVK.est.s.together}
\begin{aligned}
&\sum_{n \geq N} \left| \frac{1}{2 \pi \ii} \int_{\Gamma_n} \langle KB^{s+1} K f,f \rangle \; \dd z \right|
\leq 
\sum_{t=1}^{s+1} \binom{s+2}{t} \tau_N^{s+1} \|f\|^2		
\\
& \quad \leq  \tau_N^{s+1} \sum_{t=0}^{s+2} \binom{s+2}{t} \|f\|^2 
=
2^{s+2} \tau_N^{s+1} \|f\|^2.
\end{aligned}
\end{equation}
Since $\tau_N=o(1)$ as $N \to \infty$, we can select $N_* \in \N$ with $N_* > N_0$ such that
\begin{equation}
\tau_{N_*} \leq \frac 14.
\end{equation}
With this choice, both \eqref{tau.s0.cond} and $2^{s+2} \tau_{N_*}^{s+1} \leq 2^{-s}$ are satisfied.
\end{proof}

\begin{proof}[Proof of Theorem~\ref{thm:RB.slow}]
To justify the claim on the Riesz property, we employ Theorem~\ref{thm:RB.criteria} for the 
system of projections $\{S_0\}\cup\{P_n\}_{n > N_0}$. Recall that these are disjoint and complete; see Proposition~\ref{prop:localization}. Hence,
it suffices to verify that 
\begin{equation}\label{RB.cond.N*}
	\forall f \in \cH\, : \, \sum_{n \geq N_*} |\langle (P_n - P_n^0) f, f \rangle| < \infty,
\end{equation}
where $N_* > N_0$ is selected in Lemma~\ref{lem:KVK.s.est}. 

To this end, we employ the representation of the resolvent difference
\begin{equation}\label{res.exp.s}
(z-T)^{-1} - (z-A)^{-1} = K(z) \left(\sum_{s =0}^\infty B(z)^{s+1} \right) K(z);
\end{equation}
recall that $\|B(z)\| \leq 1/2$ for all $z \notin \Pi$ with $\Re z \geq \mu_{N_0} -  r_{N_0}$; see \eqref{KVK.1/2} in the proof of Proposition~\ref{prop:localization}.
From Lemma~\ref{lem:KVK.s.est}, we arrive at
\begin{equation}
	\sum_{s =0}^\infty\sum_{n \geq N_*} \left| \frac{1}{2 \pi \ii} \int_{\Gamma_n} \langle KB^{s+1} K f,f \rangle \; \dd z \right|
	\leq \sum_{s =0}^\infty 2^{-s} \|f\|^2 = 2 \|f\|^2.
\end{equation}
Hence 
\begin{equation}
\begin{aligned}
\sum_{n \geq N_*} |\langle (P_n - P_n^0) f, f \rangle| &= 
\sum_{n \geq N_*} \left| \frac{1}{2 \pi \ii} \int_{\Gamma_n} \Big\langle 
K \Big(\sum_{s =0}^\infty B^{s+1} \Big) K f,f 
\Big \rangle \; \dd z \right|
 < \infty,
\end{aligned}
\end{equation}
thus \eqref{RB.cond.N*} indeed holds.
\end{proof}

\section{Simpler proof if $\{\omega_j\}_{j \in \N}$ decays faster}
\label{sec:fast}

We consider the sequence $\{\omega_j\}_{j \in \N}$ in \eqref{asm:v.loc.sub} satisfying \eqref{asm:om.l1}, \eqref{asm:omega.G} and, \emph{in addition}, 
\begin{equation}\label{asm:omega.fast}
\sup_{k \in \N} \sum_{\substack{n>N_1\\ n \neq k}} \frac{r_n}{|\mu_n-\mu_k|} \sum_{j \neq n} \frac{\omega_j^2}{|\mu_n-\mu_j|} < \infty
\end{equation}
for some $N_1 \in \N$; we will further assume that $N_1>N_0$ with $N_0$ from Proposition~\ref{prop:localization}. In this case, partially inspired by \cite[Step 2, Proof of Thm.~1]{Motovilov-2017-8}, the proof of the statement in Theorem~\ref{thm:RB.slow} can be substantially simplified. 

\begin{proof}[Proof of Theorem~\ref{thm:RB.slow} with the additional assumption~\eqref{asm:omega.fast}]
In view of Theorem~\ref{thm:RB.criteria}, we justify that \eqref{RB.cond.N*} holds (the remaining steps of the proof remain unchanged). 

Instead of the resolvent expansion \eqref{res.exp.s}, it suffices to employ
\begin{equation}\label{res.dif}
(z-T)^{-1} - (z-A)^{-1} =  K(z) B(z) (I-B(z))^{-1}K(z)
\end{equation}  
and to show that, for any $f \in \cH$,
\begin{equation}\label{sum.conv}
\sum_{n>N_1}  \int_{\Gamma_n} | \langle B (I-B)^{-1} Kf, K^* f \rangle | |\dd z| < \infty.
\end{equation}
To this end, we decompose $f$ in the orthonormal basis $\{\psi_k\}_{k \in \N}$, i.e.,
\begin{equation}
	f = \sum_k f_k \psi_k,
\end{equation}
and notice that 
	\begin{equation}\label{Kf.norm}
		\|Kf\|^2 = \sum_k |\langle Kf, \psi_k \rangle|^2 = \sum_k \frac{|f_k|^2}{|z-\mu_k|} = \|K^* f\|^2.  
	\end{equation}
Since $\|B(z)\| \leq 1/2$ for $z \in \Gamma_n$ with $n > N_1 >N_0$ (see the proof of Proposition~\ref{prop:localization}),  
the sum on the l.h.s.~of \eqref{sum.conv} does not exceed (see also \eqref{KVK.HS})
\begin{equation}
\begin{aligned}
&\sum_{n >N_1}  \int_{\Gamma_n}  \frac{\|B\|}{1-\|B\|} \sum_k \frac{|f_k|^2}{|z-\mu_k|} |\dd z|
 \leq 
2\sum_{n >N_1}  \int_{\Gamma_n} \|B\| \sum_k \frac{|f_k|^2}{|z-\mu_k|} |\dd z|
\\
& \quad \leq 
2\sum_{n >N_1}  \int_{\Gamma_n} \sum_j \frac{\omega_j^2}{|z-\mu_j|} \sum_k \frac{|f_k|^2}{|z-\mu_k|} |\dd z|
\\
& \quad \leq
2 \sup_{k \in \N} 
\left(
\sum_{n >N_1}  \int_{\Gamma_n} \frac{1}{|z-\mu_k|} \sum_j \frac{\omega_j^2}{|z-\mu_j|}  |\dd z|
\right)
\|f\|^2
\\
& \quad \leq
4 \pi \sup_{k \in \N} 
\left(
\sum_{n >N_1}  \frac{r_n}{|z_n-\mu_k|} \sum_j \frac{\omega_j^2}{|z_n-\mu_j|}  
\right)
\|f\|^2,
\end{aligned}
\end{equation}
where the sequence $z_n \in \Gamma_n$, $n>N_1$, is selected such that the maximum of the integrand is achieved. Finally, arguing as in \eqref{sum.KVK.z.1} and \eqref{sum.KVK.z.3}, we infer that 
\begin{equation}
\begin{aligned}
&\sum_{n >N_1}  \frac{r_n}{|z_n-\mu_k|} \sum_j \frac{\omega_j^2}{|z_n-\mu_j|}  
\\
& \quad
\leq
C \left(
\sum_{\substack{n>N_1\\ n \neq k}} \frac{r_n}{|\mu_n-\mu_k|} \sum_{j \neq n} \frac{\omega_j^2}{|\mu_n-\mu_j|}
+ 
\sum_{j \neq k} \frac{\omega_j^2}{|\mu_k-\mu_j|}
+ \frac{\omega_k^2}{r_k}
\right).
\end{aligned}
\end{equation}
Thus, using the assumptions~\eqref{asm:omega.G} and \eqref{asm:omega.fast}, we obtain that \eqref{sum.conv} holds.
\end{proof}

\section{Sufficient conditions for \eqref{asm:omega.G} if $\mu_n=n$, $n \in \N$, and examples in this case}
\label{sec:ex}

In this section, we specialize to the case $\mu_n=n$, $n \in \N$. In Sections~\ref{ssec:om.dec} and \ref{ssec:om.gaps}, we find sufficient conditions on $\{\omega_j\}_{j\in \N}$ for \eqref{asm:omega.G}, which are expressed as a sufficiently fast decay or as an arbitrarily slow decay combined with sufficiently large gaps in the support of $\{\omega_j\}_{j\in \N}$. In Section~\ref{ssec:fin-band}, we apply these in the example of tridiagonal (or finite-band) matrices. In Section~\ref{ssec:suff.cond.fast}, we find a sufficient condition for \eqref{asm:omega.fast} in terms of the decay rate of  $\{\omega_j\}_{j\in \N}$ and compare it with the one for \eqref{asm:omega.G} found in Section~\ref{ssec:om.dec}. Finally, in Section~\ref{ssec:counter.ex}, we give an example showing that the summability condition~\eqref{asm:om.l1} does not suffice for the Riesz basis property.

\subsection{Sufficiently fast decay of $\{\omega_j\}_{j\in \N}$}
\label{ssec:om.dec}

\begin{lemma}\label{lem:sum.o(ln)}
	Let $\{t_k\}_{k \in \N} \subset \R_+$  be such that
	\begin{equation}\label{tk.o(ln)}
		\left\{\frac{t_k}{k} \right \}_{k \in \N} \in \ell^1(\N) \quad \text{and} \quad t_k \log k = o(1), \quad k \to \infty.
	\end{equation}
	Then
	\begin{equation}
		\sum_{k \neq n} \frac{t_k}{|k-n|} = o(1), \quad n \to \infty.
	\end{equation}
\end{lemma}

\begin{proof}
	First we have
	\begin{equation}
		\sum_{k \geq 2n} \frac{t_k}{|k-n|} = 	\sum_{k \geq 2n} \frac{t_k}{k-n} 
		\leq
		2 \sum_{k \geq 2n} \frac{t_k}{k} = o(1), \quad n \to \infty.
	\end{equation}
	Next, 
\begin{equation}
	\lim_{n \to \infty}\sum_{1 \leq k \leq \frac n2} \frac{t_k}{|k-n|}=
	\lim_{n \to \infty}\sum_{1 \leq k \leq \frac n2} \frac{t_k}{n-k}  \leq  \lim_{n \to \infty } \frac{2}{n} \sum_{1 \leq k \leq n} t_k 
\end{equation}
and by Stolz theorem (see e.g.~\cite[Problem~70]{Polya-1972}) we arrive at
\begin{equation}
	\lim_{n \to \infty } \frac{2}{n} \sum_{1 \leq k \leq n} t_k = 2 \lim_{n \to \infty} t_{n+1} = 0.
\end{equation}
	The estimates in the regions $\frac n2 \leq k \leq n-1$ and $n+1 \leq k \leq 2n-1$ are analogous, we give details for the latter only. Namely,

	\begin{equation}
		\begin{aligned}
			\sum_{n+1 \leq k \leq 2n-1 } \frac{t_k}{|k-n|} & =
			\sum_{n+1 \leq k \leq 2n-1 } \frac{t_k \log k}{ (k-n) \log k}
			\\
			& \leq \Big(\max_{n+1 \leq k \leq 2n-1} t_k \log k \Big) \sum_{n+1 \leq k \leq 2n-1 } \frac{1}{ (k-n) \log k}
			\\
			& 
			\leq \Big(\max_{n+1 \leq k \leq 2n-1} t_k \log k \Big) \frac 1 {\log(n+1)}\sum_{1 \leq j \leq n-1 } \frac{1}{j}
			\\
			& 
			\leq \Big(\max_{n+1 \leq k \leq 2n-1} t_k \log k \Big) \frac {1+\log (n-1)} {\log(n+1)}.
		\end{aligned}
	\end{equation}
	From $t_k \log k = o(1)$, we arrive at
	\begin{equation}
		\lim_{n \to \infty} \max_{n+1 \leq k \leq 2n-1} (t_k \log k)	= 0.
	\end{equation}
\end{proof}

\begin{remark}\label{rem:decr}
	In Lemma~\ref{lem:sum.o(ln)}, instead of $t_k \log k = o(1)$, one can assume that $\{t_k\}_{k \in \N}$ is decreasing and $t_k \to 0$ as $k \to \infty$. Indeed, since for $n>6$
	\begin{equation}
		\sum_{\sqrt n \leq k \leq n} \frac{t_k}{k} \geq t_n \sum_{\sqrt n \leq k \leq n} \frac{1}{k} \geq t_n \int_{\sqrt n+1}^n \frac{\dd x}{x} 
		= t_n  \log \frac{n}{\sqrt n+1},
	\end{equation}
	we obtain from $\{t_k/k\}_{k \in \N} \in \ell^1(\N)$ that $t_n \log n = o(1)$ as $n \to \infty$.  	Thus the assumptions of Lemma~\ref{lem:sum.o(ln)} are  satisfied.

\end{remark}

\begin{example}\label{ex:lnln}
	The sequence
	\begin{equation}\label{om.lnln}
		\omega_j =\frac{1}{(\log j)^\frac12 (\log \log j)^a }, \quad j \geq 2, \ a > \frac 12
	\end{equation}
	satisfies the conditions \eqref{asm:om.l1} and \eqref{asm:omega.G} in Theorem~\ref{thm:RB.slow} for $\mu_n=n$, $n \in \N$. 
	
	Indeed, in view of Lemma~\ref{lem:sum.o(ln)} we need to verify that $\{\omega_j^2/j\}_{j \in \N} \in \ell^1(\N)$ only, which follows from the integral test since
	\begin{equation}
		\int_{e^e}^\infty \frac{\dd x}{x \log x (\log \log x)^{2a} } = 
		\int_e^\infty \frac{\dd y}{y (\log y)^{2a} } = 
		\int_{1}^\infty \frac{\dd u}{ u^{2a} } < \infty. 
	\end{equation}
\end{example}

\subsection{An arbitrarily slow decay and sufficiently large gaps in the support of $\{\omega_j\}_{j\in \N}$}
\label{ssec:om.gaps}

Let $b$ be the sequence defined by 
\begin{equation}\label{bm.def}
	b_m:=[m^a], \quad m \in \N, \ a>1;	
\end{equation}
here $[\cdot]$ denotes the integer value. We consider the sequence $t$ such that
\begin{equation}\label{t.gap.1}
	t_j = 0 \quad \text{if} \quad j \notin \Ran(b) \quad \text{and} \quad  \lim_{m \to \infty}t_{b_m} =0;
\end{equation}
i.e.~the gaps between the non-zero elements of the sequence $t$ are determined by the sequence $b$ and we assume that the subsequence of the non-zero elements tends to $0$ (with no restrictions on the rate).
\begin{lemma}\label{lem:t.gaps}
	Let the sequences $b$ and $t$ satisfy \eqref{bm.def} and \eqref{t.gap.1}. Then
	\begin{equation}\label{Tn.dec}
		\sum_{\substack{j =1 \\ j \neq n}}^\infty \frac{t_j}{|n-j|} = \BigO \left(\frac{\log n}{n^{1-\frac 1a}} \right) + \max_{n^\frac1a-2<m < n^\frac1a+3} |t_{b_m}|, \quad n \to \infty.
	\end{equation}
\end{lemma}
\begin{proof}
	We notice first that the sum in \eqref{Tn.dec} can be rewritten as
	\begin{equation}
		\sum_{\substack{m=1 \\ b_m \neq n}}^\infty \frac{s_m}{|n-b_m|}, 
	\end{equation}
	where the sequence $s$ is defined by $s_m:= t_{b_m}$, $m \in \N$.
	For $n \in \N$, let $p = n^\frac1a$, i.e.~$p^a =n$. From the properties of the function $[\cdot]$, we have
	\begin{equation}
		\begin{aligned}
			p^a- m^a & \leq n - b_m \leq p^a- (m-1)^a, &m \leq p-2,
			\\
			(m-1)^a - p^a & \leq b_m -n \leq m^a - p^a, &m \geq p+2.
		\end{aligned}
	\end{equation}
	Thus 
	\begin{equation}\label{Sn.dec}
\begin{aligned}
\sum_{\substack{m=1 \\ b_m \neq n}}^\infty \frac{|s_m|}{|n-b_m|}
&\leq 
\sum_{1 \leq m \leq p-2} \frac{\|s\|_{\ell^\infty}}{p^a-m^a} 
+ 
\sum_{\substack{p-2<m < p+3 \\ m \neq p}} |s_m|
\\
& \qquad +
		\sum_{p+3 \leq m} \frac{\|s\|_{\ell^\infty}}{(m-1)^a - p^a}.
\end{aligned}
	\end{equation}
	The second sum on the r.h.s.~of \eqref{Sn.dec} tends to $0$ due to the second assumption in \eqref{t.gap.1}. The remaining two sums are estimated via integrals. To this end, by comparing $y$ and $y^a$ on intervals $(0,1)$ and $(1,2)$, we obtain
	\begin{equation}\label{y.ya.comp}
		\sup_{y \in (0,1) \cup (1,2) } \frac{1-y}{1-y^a} =1.
	\end{equation}
	The straightforward estimates, the substitution $x=py$ and \eqref{y.ya.comp} yield, as $n \to \infty$,
	\begin{equation}
		\begin{aligned}
			\sum_{1 \leq m \leq p-2} \frac{1}{p^a-m^a} 
			&\leq 
			\int_1^{p-1} \frac{\dd x}{p^a-x^a} 
			\leq  
			\frac{1}{p^{a-1}}\int_0^{1-\frac 1p} \frac{1-y}{1-y}\frac{\dd y}{1-y^a} 
			\\
			& \leq 
			\frac{1}{p^{a-1}}
			\int_0^{1-\frac 1p} \frac{\dd y}{1-y} 
			= \frac{\log p}{p^{a-1}} 
			= \BigO \left(\frac{\log n}{n^{1-\frac 1a}} \right).
		\end{aligned}
	\end{equation}
	Similarly, as $n \to \infty$,
	\begin{equation}
		\begin{aligned}
			\sum_{p+3 \leq m \leq 2p+1} \frac{1}{(m-1)^a - p^a}
			&\leq 
			\int_{p+2}^{2p+1} \frac{\dd x}{(x-1)^a - p^a} 
			=  
			\frac{1}{p^{a-1}}\int_{1+\frac1p}^{2}\frac{\dd y}{y^a-1} 
			\\
			& \leq 
			\frac{1}{p^{a-1}}
			\int_{1+\frac 1p}^2 \frac{\dd y}{y-1} 
			= \BigO \left(\frac{\log n}{n^{1-\frac 1a}} \right), 
			\\
			\sum_{2p +2\leq m} \frac{1}{(m-1)^a - p^a}
			&\leq 
			\int_{2p}^{\infty} \frac{\dd x}{x^a - p^a} 
			\leq \frac{2^a}{2^a- 1} \int_{2p}^{\infty} \frac{\dd x}{x^a} = \BigO(n^{1-a})
		\end{aligned}
	\end{equation}
	and so the claim  \eqref{Tn.dec} follows.
\end{proof}
\begin{example}\label{ex:t.gaps}
	Let $\{\omega_j\}_{j \in \N} \subset [0,\infty)$ be such that the assumption of Lemma~\ref{lem:t.gaps} hold for $t_j = \omega_j^2$, $j \in \N$, and some $a>1$ in \eqref{bm.def}. Then, with such $\{\omega_j\}_{j \in \N}$, the conditions \eqref{asm:om.l1} and \eqref{asm:omega.G} in Theorem~\ref{thm:RB.slow} for $\mu_n=n$, $n \in \N$, are satisfied. 
	
	Indeed, in view of Lemma~\ref{lem:t.gaps}, we need to verify only that $\{\omega_j^2/j\}_{j \in \N}$ is summable. To this end, notice that
	\begin{equation}
		\sum_{j=1}^\infty \frac{\omega_j^2}{j} = \sum_{m=1}^\infty \frac{\omega_{b_m}^2}{b_m}
		\leq \|\{\omega_j\}_{j \in \N}\|_{\ell^\infty} \left(1+ \sum_{m=2}^\infty \frac 1{(m-1)^a} \right) < \infty.
	\end{equation}
	
\end{example}

\subsection{Application in finite-band matrices}
\label{ssec:fin-band}

The sufficient conditions from Subsections~\ref{ssec:om.dec} and \ref{ssec:om.gaps} can be used to analyze tridiagonal (or finite-band) perturbations of diagonal matrices in $\ell^2(\N)$.

\begin{example}
	Let $\cH=\ell^2(\N)$ and 
	\begin{equation}
		A =
		\begin{pmatrix}
			1 & 0 & 0 & 0 & . 
			\\
			0 & 2 & 0 & 0 & . 
			\\
			0 & 0 & 3 & 0 & . 
			\\
			0 & 0 & 0 & 4 & . 
			\\
			. & . & . & . & . 
		\end{pmatrix},
		\qquad
		V =
		\begin{pmatrix}
			b_1^{(0)} & b_1^{(1)} & 0 & 0 & . 
			\\
			b_1^{(-1)} & b_2^{(0)} & b_2^{(1)} & 0 & . 
			\\
			0 & b_2^{(-1)} & b_3^{(0)} & b_3^{(1)} & . 
			\\
			0 & 0 & b_3^{(-1)} & b_4^{(0)} & . 
			\\
			. & . & . & . & . 
		\end{pmatrix},
	\end{equation}
	where, for $j\in \{-1,0,1\}$, the sequences $\{b_k^{(j)}\}_{k \in \N}$ satisfy
	\begin{equation}\label{bjk.est}
		|b_{k-1}^{(j)}|+|b_k^{(j)}| = \BigO(\omega_k^2), \quad k \to \infty,
	\end{equation}
	with a non-negative sequence $\{\omega_k\}_{k \in \N}$. In this case, the eigenvectors of $A$ read $\psi_k = e_k$, $k \in \N$, where $\{e_k\}_{k \in \N}$ is the standard basis of $\ell^2(\N)$. From the condition \eqref{bjk.est}, we obtain that there is $C>0$ such that for all sufficiently large $m,n \in \N$
	\begin{equation}
		|\langle V e_m, e_n \rangle| \leq C \min \{\omega_m^2, \omega_n^2\} 
		\leq C \omega_m \omega_n.
	\end{equation}
	Suppose further that
	\begin{equation}
		\omega_j = \omega_j^{(1)} + \omega_j^{(2)}, \quad j \in \N,
	\end{equation}
	where with some $a >1/2$
	\begin{equation}
		\omega_j^{(1)} = \BigO \Big( \frac{1}{\log j(\log \log j)^a} \Big), \qquad j \to \infty,	
	\end{equation}
	cf.~Example~\ref{ex:lnln}, and $\{\omega_j^{(2)}\}_{j \in \N}$ be such that $t_j = (\omega_j^{(2)})^2$, $j \in \N$, satisfies the assumptions of Lemma~\ref{lem:t.gaps} with some $a>1$, cf.~Example~\ref{ex:t.gaps}. 
	
	Notice that both $\{\omega_j^{(1)}\}_{j \in \N}$ and $\{\omega_j^{(2)}\}_{j \in \N}$ satisfy the conditions \eqref{asm:om.l1} and \eqref{asm:omega.G} and, by a simple inequality, also $\{\omega_j\}_{j \in \N}$ satisfies both conditions. Hence the assumptions of Theorem~\ref{thm:RB.slow} are satisfied and its claim holds for the operator $T=A+V$. 
	
	We remark that if $\omega_j^{(1)}=0$, $j \in \N$,  the large gaps in $\{\omega_j^{(2)}\}_{j \in \N}$ allow for a direct analysis of $A+V$ since it has a block diagonal structure. However, such a reduction is not possible with a non-trivial $\{\omega_j^{(1)}\}_{j \in \N}$.
	
	Finally, this example can be naturally extended to a finite band matrices with conditions on the perturbations analogous to \eqref{bjk.est}, cf.~also \cite[Sec.~5.2]{Mityagin-2019-139}.
	
\end{example}

\subsection{Sufficiently fast decay of $\{\omega_j\}_{j\in \N}$ for \eqref{asm:omega.fast}}
\label{ssec:suff.cond.fast}

\begin{lemma}\label{lem:log.est}
	Let $\gamma>1$. Then
	\begin{equation}\label{sum.ln.gamma}
		\sum_{2<k \neq n} \frac{1}{(\log k)^\gamma|n-k|} = \BigO \left(\frac 1{(\log n)^{\gamma-1}} \right), \qquad n \to \infty.
	\end{equation}
\end{lemma}
\begin{proof}
	We split the estimate of the sum in \eqref{sum.ln.gamma} to several parts. First, as $n \to \infty$,
	\begin{equation}
		\sum_{2< k < \sqrt n} \frac{1}{(\log k)^\gamma|n-k|} \leq \frac{\sqrt n}{(\log 3)^\gamma (n-\sqrt n)} = \BigO(n^{-\frac 12}). 
	\end{equation}
	Next,
	\begin{equation}
		\begin{aligned}
			\sum_{\sqrt n \leq k \leq n-1 } \frac{1}{(\log k)^\gamma|n-k|} 
			&\leq \frac{1}{(\log \sqrt n)^\gamma} \sum_{1 \leq j \leq n} \frac{1}{j} = \BigO \left(\frac{1}{(\log n)^{\gamma-1}} \right). 
			\\
			\sum_{n+1 \leq k \leq 2n} \frac{1}{(\log k)^\gamma|n-k|} 
			&\leq \frac 1 {(\log n)^\gamma} \sum_{1 \leq j \leq n \frac{1}{j}}
			=\BigO \left(\frac{1}{(\log n)^{\gamma-1}} \right).  
		\end{aligned}
	\end{equation}
	Finally,
	\begin{align}
		\sum_{k \geq 2n+1} \frac{1}{(\log k)^\gamma|n-k|} & \leq 
		\int_{2n}^{\infty} \frac{\dd x}{{(\log x)^\gamma}(x-n)}
		\leq 
		2 \int_{2n}^{\infty} \frac{\dd x}{{(\log x)^\gamma}x}
		\\
		& \leq 
		2 \int_{\log (2n)}^{\infty} \frac{\dd y}{{y^\gamma}}
		= 
		\BigO \left(\frac{1}{(\log n)^{\gamma-1}} \right).
		\qedhere
	\end{align}
\end{proof}
\begin{example}\label{ex:ln.fast}
	The sequence
	\begin{equation}
		\omega_j = \frac{1}{ (\log j)^a}, \qquad j \geq 2, \ a > 1
	\end{equation}
	satisfies the assumptions of Theorem~\ref{thm:RB.slow} and also \eqref{asm:omega.fast} for $\mu_n=n$, $n \in \N$. 
	
	The former can be checked using Lemma~\ref{lem:sum.o(ln)} and the integral test (see also Example \ref{ex:lnln}). The latter follows by Lemma~\ref{lem:log.est}, indeed,
	\begin{equation}
		\begin{aligned}
			\sum_{\substack{n>2\\ k \neq n}} \frac{1}{|n-k|} \sum_{\substack{j>2\\ j \neq n}} \frac{\omega_j^2}{|n-j|}
			&\leq \sum_{\substack{n>2\\ k \neq n}} \frac 1{|n-k|}\frac{C}{(\log n)^{2a-1}} = \BigO \left(\frac{1}{(\log k)^{2(a-1)}} \right) 
		\end{aligned}
	\end{equation}
	as $k \to \infty$.

	Notice that \eqref{asm:omega.fast} is not satisfied if $a=1$. In detail, for all sufficiently large $n$,
	\begin{equation}
		\begin{aligned}
			\sum_{\substack{j>2\\ j \neq n}} \frac{\omega_j^2}{|n-j|}
			&\geq 
			\frac 1{(\log 2n)^2}\sum_{n+1 \leq j \leq 2n} \frac{1}{j-n}
			= 
			\frac 1{(\log 2n)^2}\sum_{1 \leq k \leq n} \frac{1}{k}
			\\
			&\geq \frac{\log n}{(\log 2n)^2} \geq \frac 12 \frac 1{\log n},
		\end{aligned}	
	\end{equation}
	hence,
	\begin{equation}
		\begin{aligned}
			\sum_{\substack{n>2\\ k \neq n}} \frac{1}{|n-k|} \sum_{\substack{j>2\\ j \neq n}} \frac{\omega_j^2}{|n-j|}
			&\geq \frac 12 \sum_{k \geq 2 n} \frac 1{k-n}\frac{1}{\log n} = \infty.
		\end{aligned}
	\end{equation}
Notice also that a slower decay is possible in Example~\ref{ex:lnln} where Theorem~\ref{thm:RB.slow} applies (\emph{without} the additional assumption \eqref{asm:omega.fast}).

\end{example}

\subsection{Summability condition~\eqref{asm:om.l1}}
\label{ssec:counter.ex}

The following example shows that the summability condition~\eqref{asm:om.l1} alone does not guarantee the Riesz property of the spectral projections of $T$.

\begin{example}\label{ex:counter}
We adapt the example from \cite[Sec.~6.3]{Adduci-2012-10}. Consider $\cH = \ell^2(\N)$, its standard basis $\{e_k\}_{k \in \N}$ and define $A e_k := k e_k$, $k\in \N$. We define the perturbation
\begin{equation}
	V e_{2k-1} := - \frac 12 t_{2k-1} e_{2k}, \quad  V e_{2k} := \frac 12  t_{2k} e_{2k-1}, \quad k \in \N,
\end{equation}
where $\{t_k\}_{k \in \N} \subset (0,1)$ with 
\begin{equation}\label{tk.def}
t_{2k-1} = t_{2k}=s_k, \quad k \in \N,
\end{equation}
is chosen below.

We have
\begin{equation}
	|\langle V e_m, e_n \rangle| \leq \min \{ \|Ve_m\|, \|V^*e_n\|\} \leq \|Ve_m\|^\frac 12 \|V^*e_n\|^\frac 12 \leq \frac 12 t_m^\frac 12 t_n^\frac 12,
\end{equation}
so we define $\{\omega_k\}_{k \in \N}$ by
\begin{equation}
\omega_k^2 = \frac 12 t_k, \quad k \in \N.	
\end{equation}

We select a bounded $\{s_k\}_{k \in \N}$ as
\begin{equation}\label{sk.def}
s_k = 
\begin{cases}
0, & k \neq m^2 \text{ for any } m \in \N,	
\\
\sqrt{1-\frac1k}, & k = m^2 \text{ for some } m \in \N.
\end{cases}
\end{equation}
With this choice
\begin{equation}
\begin{aligned}
\sum_{k =1}^\infty \frac{\omega_k^2}{k} &= \frac 12	\sum_{k =1}^\infty \frac{s_k}{2k-1} + \frac 12 \sum_{k =1}^\infty \frac{s_k}{2k}
\\
&\leq \frac 12\|\{s_k\}_{k \in \N}\|_{\ell^\infty} \left(
\sum_{m =1}^\infty \frac{1}{2m^2-1} + \sum_{m =1}^\infty \frac{1}{2m^2}
\right)
< \infty.
\end{aligned}
\end{equation}

Next we show directly that the eigensystem of $T=A+V$ does not contain a basis. Since the perturbation $V$ is block-diagonal, it suffices to analyze the 2-dimensional blocks corresponding to $\lspan\{e_{2k-1},e_{2k}\}$,
\begin{equation}
	T_k=A_k + V_k =
	\begin{pmatrix}
		2k- \frac 12 & 0 \\
		0 & 2k- \frac 12
	\end{pmatrix}
	+
	\frac 12
	\begin{pmatrix}
		-1  & s_k \\
		-s_k & 1
	\end{pmatrix};
\end{equation}
notice that for the adjoint operator we have the blocks
\begin{equation}
	T_k^*=A_k + V_k^* =
	\begin{pmatrix}
		2k- \frac 12 & 0 \\
		0 & 2k- \frac 12
	\end{pmatrix}
	+
	\frac 12
	\begin{pmatrix}
		-1  & - s_k \\
		s_k & 1
	\end{pmatrix}.
\end{equation}

Eigenvalues and eigenvectors of $T_k$ and $T_k^*$ can be calculated explicitly, namely,
\begin{equation}
	\begin{aligned}
		T_k g_k^\pm &= (2k- \tfrac12 \pm \tfrac 12 \varsigma_k) g_k^\pm, 
		\\
		T_k^* (g_k^*)^\pm &= (2k- \tfrac12 \pm \tfrac 12 \varsigma_k) (g_k^*)^\pm, 
		\\ 
		g_k^\pm &= \begin{pmatrix}
			1 \\ G_k^{\pm 1}
		\end{pmatrix},
		\quad
		(g_k^*)^\pm = \begin{pmatrix}
			1 \\ - G_k^{\pm 1}
		\end{pmatrix},
		\\
		G_k &= \left(\frac{1+\varsigma_k}{1-\varsigma_k} \right)^\frac 12, \quad \varsigma_k = \sqrt{1-s_k^2}. 
	\end{aligned}
\end{equation}
The norms of the one-dimensional spectral projections $P_k^\pm$ of $T$ related to eigenvalues $2k - \tfrac 12 \pm \tfrac 12 \varsigma_k$ read
\begin{equation}
\|P_k^\pm\| = \frac{\| \langle \cdot, (g_k^*)^\pm \rangle g_k^\pm \|}{|\langle g_k^\pm, (g_k^*)^\pm \rangle  |} = 
\frac{1+G_k^2}{1-G_k^2}
=
\frac{1}{\varsigma_k} = \frac{1}{\sqrt{1 - s_k^2}}.
\end{equation}
Hence, with our choice of the sequence $\{s_k\}_{k \in \N}$ and $\{t_k\}_{k \in \N}$ in \eqref{sk.def} and \eqref{tk.def}, $\varsigma_{m^2} = 1/m$ and we arrive at 
\begin{equation}
\|P_{m^2}^\pm\| = m, \qquad m \in \N.
\end{equation}
\end{example}

{
	\bibliographystyle{halpha}
	\bibliography{C:/Data/00Synchronized/references}
}

\end{document}